\documentclass{amsart}
\usepackage{amssymb,euscript,amsmath, mathrsfs}
\usepackage[dvips]{graphicx}
\usepackage[dvips]{color}
\usepackage{epsfig}

\newcounter{ENUM}
\newcommand{\itm}{\item}
\newenvironment{ilist}{\renewcommand{\theENUM}{\roman{ENUM}}\renewcommand{\itm}{\addtocounter{ENUM}{1}\item[(\theENUM)]}\begin{itemize}\setcounter{ENUM}{0}}{\end{itemize}}
\newenvironment{alist}[1][0]{\renewcommand{\theENUM}{\alph{ENUM}}\renewcommand{\itm}{\addtocounter{ENUM}{1}\item[\theENUM)]}\begin{itemize}\setcounter{ENUM}{#1}}{\end{itemize}}

\newcommand{\nHat}[1]{\widehat{#1}}                        
\newcommand{\nBar}[1]{\overline{#1}}

\def\S{{\mathfrak S}}
\def\x{{\mathbf x}}
\def\s{{\mathbf s}}
\def\1{{\mathbf 1}}

\def\Z{{\mathbb Z}}
\def\N{{\mathbb N}}

\def\R{{\mathbb R}}

\def\L{{\mathcal L}}

\def\m{{\mathfrak m}}

\def\sgn{\operatorname{sign}}

\def\conv{\mathrm{conv}}
\def\vol{\mathrm{Vol}}

\newtheorem{thm}{Theorem}[section]
\newtheorem{prop}[thm]{Proposition}
\newtheorem{lem}[thm]{Lemma}
\newtheorem{cor}[thm]{Corollary}
\newtheorem{conj}[thm]{Conjecture}

\theoremstyle{definition}
\newtheorem{defn}[thm]{Definition}
\newtheorem{ques}[thm]{Question}
\newtheorem{ex}[thm]{Example}

\theoremstyle{remark}

\newtheorem{rem}[thm]{Remark}

\numberwithin{equation}{section}

\address{Room 2-333, 77 Massachusetts Avenue, Department of Mathematics, Massachusetts Institute of Technology, Cambridge, MA 02139, USA}
\email{fuliu@math.mit.edu}

\begin{document}
\title{Ehrhart polynomials of lattice-face polytopes}
\author{Fu Liu}
\begin{abstract}
There is a simple formula for the Ehrhart polynomial of a cyclic
polytope. The purpose of this paper is to show that the same formula
holds for a more general class of polytopes, lattice-face polytopes.
We develop a way of decomposing any $d$-dimensional simplex in
general position into $d!$ signed sets, each of which corresponds to
a permutation in the symmetric group $\S_d,$ and reduce the problem
of counting lattice points in a polytope in general position to that
of counting lattice points in these special signed sets. Applying
this decomposition to a lattice-face simplex, we obtain signed sets
with special properties that allow us to count the number of lattice
points inside them. We are thus able to conclude the desired formula
for the Ehrhart polynomials of lattice-face polytopes.
\end{abstract}

\maketitle

\section{Introduction}
A $d$-dimensional {\it lattice} $\Z^d = \{\x = (x_1, \dots, x_d) \ |
\ \forall x_i \in \Z\}$ is the collection of all points with integer
coordinates in $\R^d.$ Any point in a lattice is called a {\it
lattice point}.

A {\it convex polytope} is a convex hull of a finite set of points.
We often omit convex and just say polytope. For any polytope $P$ and
some positive integer $m \in \N,$ we use $i(m,P)$ to denote the
number of lattice points in $m P,$ where $m P = \{ m x | x \in P \}$
is the {\it $m$th dilated polytope} of $P.$

An {\it integral} or {\it lattice} polytope is a convex polytope
whose vertices are all lattice points. Eug\`{e}ne Ehrhart
\cite{Ehrhart} showed that for any $d$-dimensional integral
polytope, $i(P,m)$ is a polynomial in $m$ of degree $d.$ Thus, we
call $i(P,m)$ the Ehrhart polynomial of $P$ when $P$ is an integral
polytope. Although Ehrhart's theory was developed in the 1960's, we
still do not know much about the coefficients of Ehrhart polynomials
for general polytopes except that the leading, second and last
coefficients of $i(P,m)$ are the normalized volume of $P$, one half
of the normalized volume of the boundary of $P$ and $1,$
respectively.

In \cite{cyclic}, the author showed that for any $d$-dimensional
cyclic polytope $P$, we have that
\begin{equation}\label{ques}
i(P, m) =  \vol(mP) + i(\pi(P), m) = \sum_{k=0}^d
\vol_k(\pi^{(d-k)}(P)) m^k,
\end{equation}
where $\pi^{(k)}: \R^d \to \R^{d-k}$ is the map which ignores the
last $k$ coordinates of a point, and asked whether there are other
integral polytopes that have the the same form of Ehrhart
polynomials.

In this paper, we define a new family of integral polytopes, {\it
lattice-face} polytopes, and show (Theorem \ref{main1}) that their
Ehrhart polynomials are in the form of (\ref{ques}).

The main method of  \cite{cyclic} is a decomposition of an arbitrary
$d$-dimensional simplex cyclic polytope into $d!$ signed sets, each
of which corresponds to a permutation in the symmetric group $\S_d$
and has the same sign as the corresponding permutation. However, for
general polytopes, such a decomposition does not work.

In this paper, we develop a way of decomposing any $d$-dimensional
simplex in {\it general position} into $d!$ signed sets, where the
sign of each set is not necessarily the same as the corresponding
permutation. Applying the new decomposition to a lattice-face
simplex, we are able to show (Theorem \ref{main2}) that the number
of lattice points is given by in terms of a formula (\ref{gsigma})
involving Bernoulli polynomials, signs of permutations, and
determinants, and then to analyze this formula further to derive the
theorem. Theorem \ref{main2}, together with some simple observations
in section 2 and 3, implies Theorem \ref{main1}.

\section{Preliminaries}

We first give some definitions and notation, most of which follows
\cite{cyclic}.

All polytopes we will consider are full-dimensional, so for any
convex polytope $P,$ we use $d$ to denote both the dimension of the
ambient space $\R^d$ and the dimension of $P.$ We call a
$d$-dimensional polytope a $d$-polytope. Also, We use $\partial P$
and $I(P)$ to denote the boundary and the interior of $P,$
respectively.



For any set $S,$ we use $\conv(S)$ to denote the convex hull of all
of points in $S.$

 Recall that the projection $\pi: \R^d
\to \R^{d-1}$ is the map that forgets the last coordinate. For any
set $S \subset \R^d$ and any point $y \in \R^{d-1},$ let $\rho(y, S)
= \pi^{-1}(y) \cap S$ be the intersection of $S$ with the inverse
image of $y$ under $\pi.$ Let $p(y, S)$ and $n(y, S)$ be the point
in $\rho(y,S)$ with the largest and smallest last coordinate,
respectively. If $\rho(y,S)$ is the empty set, i.e., $y \not\in
\pi(S),$ then let $p(y, S)$ and $n(y,S)$ be empty sets as well.
Clearly, if $S$ is a $d$-polytope, $p(y, S)$ and $n(y, S)$ are on
the boundary of $S.$ Also, we let $\rho^+(y,S) = \rho(y,S)\setminus
n(y,S),$ and for any $T \subset \R^{d-1},$ $\rho^+(T,S) = \cup_{y
\in T} \rho^+(y,S).$

\begin{defn}\label{dfbdry}
Define $PB(P) = \bigcup_{y \in \pi(P)} p(y,P)$ to be the {\it
positive boundary} of $P;$ $NB(P) = \cup_{y \in \pi(P)} n(y,P)$ to
be the {\it negative boundary} of $P$ and $\Omega(P) = P \setminus
NB(P) = \rho^+(\pi(P), P) = \cup_{y \in \pi(P)} \rho^+(y,P)$ to be
the {\it nonnegative part} of $P.$
\end{defn}
\begin{defn}
For any facet $F$ of $P,$ if $F$ has an interior point in the
positive boundary of $P,$  then we call $F$ a {\it positive facet}
of $P$ and define the sign of $F$ as $+1: \sgn(F) = + 1.$ Similarly,
we can define the {\it negative facets} of $P$ with associated sign
$-1.$  For the facets that are neither positive nor negative, we
call them {\it neutral facets} and define the sign as $0.$
\end{defn}
It's easy to see that $F \subset PB(P)$ if $F$ is a positive facet
and $F \subset NB(P)$ if $F$ is a negative facet.

We write $P = \bigsqcup_{i=1}^k P_i$ if $P = \bigcup_{i=1}^k P_i$
and for any $i \neq j$, $P_i \cap P_j$ is contained in their
boundaries. If $F_1, F_2, \dots, F_{\ell}$ are all the positive
facets of $P$ and $F_{\ell +1}, \dots, F_k$ are all the negative
facets of $P,$ then $$\pi(P) = \bigsqcup_{i=1}^{\ell} \pi(F_i) =
\bigsqcup_{i=\ell+1}^k \pi(F_i).$$

 Because the usual set union
and set minus operation do not count the number of occurrences of an
element, which is important in our paper, from now on we will
consider any polytopes or sets as {\it multisets} which allow {\it
negative multiplicities.} In other words, we consider any element of
a multiset as a pair $(\x, m),$ where $m$ is the multiplicity of
element $\x.$ Then for any multisets $M_1, M_2$ and any integers
$m,n$ and $i,$ we define the following operators:
\begin{alist}
\itm Scalar product: $i M_1 = i \cdot M_1 = \{ (\x, i m) \ | \ (\x,
m) \in M_1\}.$ \itm Addition: $M_1 \oplus M_2 = \{ (\x, m + n) \ | \
(\x, m) \in M_1, (\x, n) \in M_2 \}.$ \itm Subtraction: $M_1 \ominus
M_2 = M_1 \oplus ((-1) \cdot M_2).$
\end{alist}

It's clear that the following holds:
\begin{lem} \
For any polytope $P \subset \R^d,$  $\forall R_1, \dots, R_k \subset
\R^{d-1},$ $\forall i_1,\dots, i_k \in \Z:$
$$\rho^+\left( \bigoplus_{j=1}^k i_j R_j, \ P \right) =
\bigoplus_{j=1}^k i_j\rho^+( R_j,\ P). $$
%
\end{lem}

\begin{defn}
We say a set $S$ has {\it weight} $w,$ if each of its elements has
multiplicity either  $0$ or $w.$ And $S$ is a {\it signed} set if it
has weight $1$ or $-1.$
\end{defn}

Let $P$ be a convex polytope. For any $y$ an interior point of
$\pi(P),$ since $\pi$ is a continuous open map, the inverse image of
$y$ contains an interior point of $P.$ Thus $\pi^{-1}(y)$ intersects
the boundary of $P$ exactly twice. For any $y$ a boundary point of
$\pi(P),$ again because $\pi$ is an open map, we have that $\rho(y,
P) \subset \partial P,$ so $\rho(y,P) = \pi^{-1}(y) \cap \partial P$
is either one point or a line segment. We will only consider
polytopes $P$ where $\rho(y,P)$ always has only one point for a
boundary point $y.$

\begin{lem}\label{lemone}
If a polytope $P$ satisfies:
\begin{equation}\label{one}
|\rho(y,P)| = 1, \forall y \in \partial \pi(P),
\end{equation}
then $P$ has the following properties:
\begin{ilist}
\itm For any $y \in I(\pi(P)),$ $\pi^{-1}(y) \cap \partial P =
\{p(y,P), n(y,P)\}.$

\itm For any $y \in \partial \pi(P),$ $\pi^{-1}(y) \cap \partial P =
\rho(y,P) = p(y, P) = n(y, P),$ so $\rho^+(y,P) = \emptyset.$

\itm Let $R$ be a region containing $I(\pi(P)).$ Then $$\Omega(P) =
\rho^+(R, P) = \bigoplus_{y \in R} \rho^+(y, P).$$

\itm If $P = \bigsqcup_{i=1}^k P_i,$ where the $P_i$'s all satisfy
(\ref{one}), then $\Omega(P) = \bigoplus_{i=1}^k \Omega(P_i).$

 \itm The set of facets of $P$ are
partitioned into the set of positive facets and the set of negative
facets, i.e., there is no neutral facets.

\itm $\pi$ gives a bijection between $PB(P) \cap NB(P)$ and
$\partial \pi(P).$

\end{ilist}
\end{lem}

The proof of this lemma is straightforward, so we won't include it
here.

The main purpose of this paper is to discuss the number of lattice
points in a polytope. Therefore, for simplicity, for any set $S \in
\R^d,$ we denote by $\L(S) = S \cap \Z^d$ the set of lattice points
in $S.$ It's not hard to see that $\L$ commutes with some of the
operations we defined earlier, e.g. $\rho, \rho^+, \Omega.$

\section{Lattice-face polytopes}
A {\it $d$-simplex} is a polytope given as the convex hull of $d+1$
affinely independent points in $\R^d.$

\begin{defn}\label{dflf1} We define {\it lattice-face} polytopes recursively.
We call a one dimensional polytope a {\it lattice-face} polytope if
it is integral.

For $d \ge 2,$ we call a $d$-dimensional polytope $P$ with vertex
set $V$ a {\it lattice-face} polytope if for any $d$-subset $U
\subset V,$
\begin{alist}
\itm $\pi(\conv(U))$ is a lattice-face polytope, and

\itm $\pi(\L(H_{U})) = \Z^{d-1},$ where $H_U$ is the affine space
spanned by $U.$ In other words, after dropping the last coordinate
of the lattice of $H_U,$ we get the $(d-1)$-dimensional lattice.
\end{alist}
\end{defn}

To understand the definition, let's look at examples of
$2$-polytopes.
\begin{ex} Let $P_1$ be the polytope with vertices $v_1 = (0,0),
v_2=(2,0)$ and $v_3=(2,1).$ Clearly, for any $2$-subset $U,$
condition $a)$ is always satisfied.  When $U = \{v_1, v_2\},$ $H_U$
is $\{(x,0) \ | \ x \in \R\}.$ So $\pi(\L(H_{U})) = \Z,$ i.e., $b)$
holds. When $U = \{v_1, v_3\},$ $H_U$ is $\{ (x,y) \ | \ x = 2y\}.$
Then $\L(H_U) = \{ (2y, y) \ | \ y \in \Z \} \Rightarrow
\pi(\L(H_{U})) = 2\Z \neq \Z.$ When $U = \{v_2, v_3\},$ $H_U$ is $\{
(2,y) \ | y \in \R \ \}.$ Then $\pi(\L(H_{U})) = \{ 2\}\neq \Z.$
Therefore, $P_1$ is not a lattice-face polytope.

Let $P_2$ be the polytope with vertices $(0,0), (1,1)$ and $(2,0).$
One can check that $P_2$ is a lattice-face polytope.

\end{ex}

The following lemma gives some properties of a lattice-face
polytope.

\begin{lem}\label{plf}
 Let $P$ be a lattice-face $d$-polytope with vertex set
$V,$ then we have:
\begin{ilist}
\itm $\pi(P)$ is a lattice-face $(d-1)$-polytope.

\itm $m P$ is a lattice-face $d$-polytope, for any positive integer
$m.$

\itm $\pi$ induces a bijection between $\L(NB(P))$ (or $\L(PB(P))$)
and $\L(\pi(P)).$

\itm $\pi(\L(P)) = \L(\pi(P)).$

\itm Any $d$-subset $U$ of $V$ forms a $(d-1)$-simplex. Thus
$\pi(\conv(U))$ is a $(d-1)$-simplex.

\itm Let $H$ be the affine space spanned by some $d$-subset of $V.$
Then for any lattice point $y \in \Z^{d-1},$ we have that $\rho(y,
H)$ is a lattice point.

\itm $P$ is an integral polytope.
\end{ilist}
\end{lem}
\begin{proof}
$\rm{(i),(ii),(v)}$ and $\rm{(vi)}$ can be checked directly from the
conditions $a)$ and $b)$ of the definition. $\rm{(iii)}$ and
$\rm{(iv)}$ both follow from $\rm{(vi)}.$ We prove $\rm{(vii)}$ by
induction on $d.$

Any $1$-dimensional lattice-face polytope is integral by definition.

For $d \ge 2,$ suppose any $(d-1)$ dimensional lattice-face polytope
is an integral polytope. Let $P$ be a $d$ dimensional lattice-face
polytope with vertex set $V.$ For any vertex $v_0 \in V,$ let $U$ be
a subset of $V$ that contains $v_0.$ Let $U = \{v_0, v_1, \dots,
v_{d-1}\}.$ We know that $P' = \pi(\conv(U))$ is a lattice-face
$(d-1)$-simplex with vertices  $\{ \pi(v_0), \dots, \pi(v_{d-1})\}.$
Thus, by the induction hypothesis, $P'$ is an integral polytope. In
particular, $\pi(v_0)$ is a lattice point. Therefore, $v_0 =
\rho(\pi(v_0), H_U)$ is a lattice point.
\end{proof}

\begin{rem}
One sees that condition {\rm b)} in the definition of lattice-face
polytopes is equivalent to {\rm (vi)}.
\end{rem}

The main theorem of this paper is to describe all of the
coefficients of the Ehrhart polynomial of a lattice-face polytope.

\begin{thm}\label{main1}
 Let P be a lattice-face $d$-polytope, then
\begin{equation}
i(P, m) =  \vol(mP) + i(\pi(P), m) = \sum_{k=0}^d
\vol_k(\pi^{(d-k)}(P)) m^k.
\end{equation}

\end{thm}

However, by Lemma \ref{plf}/\rm{(ii),(iii)}, we have that
$$i(P,m) = |\L(\Omega(mP))| + i(\pi(P),m).$$

Therefore, by Lemma \ref{plf}/{\rm (i)}, to prove Theorem
\ref{main1} it is sufficient to prove the following theorem:

\begin{thm}\label{main2} For any $P$ a lattice-face polytope,
$$|\L(\Omega(P))| =   \vol(P).$$

\end{thm}


\begin{rem}\label{dflf2} We have an alternative definition of lattice-face
polytopes, which is equivalent to Definition \ref{dflf1}. Indeed, a
$d$-polytope on a vertex set $V$ is a lattice-face polytope if and
only if  for all $k$ with $0 \le k \le d-1,$
\begin{equation}\label{dfcon}
\mbox{for any $(k+1)$-subset $U \subset V,$ $\pi^{d-k}(\L(H_{U})) =
\Z^k,$} \end{equation}
 where $H_U$ is the affine space
spanned by $U.$ In other words, after dropping the last $d-k$
coordinates of the lattice of $H_U,$ we get the $k$-dimensional
lattice.
\end{rem}

\section{A signed decomposition of the nonnegative part of a simplex in general position}
The volume of a polytope is not very hard to characterize. So our
main problem is to find a way to describe the number of lattice
points in the nonnegative part of a lattice-face polytope. We are
going to do this via a signed decomposition.

\subsection{Polytopes in general position}
For the decomposition, we will work with a more general type of
polytope (which contains the family of lattice-face polytopes).

\begin{defn} We say that a $d$-polytope $P$ with vertex set $V$ is in {\it
general position} if for any $k: 0 \le k \le d-1,$ and any
$(k+1)$-subset $U \subset V,$ $\pi^{d-k}(\conv(U))$ is a
$k$-simplex, where $\conv(U)$ is the convex hull of all of points in
$U.$
\end{defn}

By the alternative definition of lattice polytopes in Remark
\ref{dflf2}, it's easy to see that a lattice-face polytope is a
polytope in general position. Therefore, the following discussion
can be applied to lattice-face polytopes.

The following lemma states some properties of a polytope in general
position. The proof is omitted.

\begin{lem}\label{pgen}
Given a $d$-polytope $P$ in general position with vertex set $V,$
then
\begin{ilist}
\itm $P$ satisfies (\ref{one}).

\itm $\pi(P)$ is a $(d-1)$-polytope in general position.

\itm For any nonempty subset $U$ of $V,$ let $Q = \conv(U).$ If $U$
is has dimension $k (0 \le k \le d),$ then $\pi^{d-k}(Q)$ is a
$k$-polytope in general position. In particular, for any facet $F$
of $P,$ $\pi(F)$ is a $(d-1)$-polytope in general position.

\itm For any triangulation of $P = \bigsqcup_{i=1}^k P_i$ without
introducing new vertices, $\Omega(P) = \bigoplus_{i=1}^k
\Omega(P_i).$ Thus, $\L(\Omega(P)) = \bigoplus_{i=1}^k
\L(\Omega(P_i)).$

\itm If $F_1, F_2, \dots, F_{\ell}$ are all the positive facets of
$P$ and $F_{\ell +1}, \dots, F_k$ are all the negative facets of
$P,$ then $\Omega(\pi(P)) = \bigoplus_{i=1}^{\ell} \Omega(\pi(F_i))
= \bigoplus_{i=\ell+1}^k \Omega(\pi(F_i)).$

\itm For any hyperplane $H$ determined by one facet of $P$ and any
$y \in \R^{d-1},$ $\rho(y, H)$ is one point.

\itm For any $k: 0 \le k \le d-1,$ any $(k+1)$-subset $U$ of $V,$
any $y_1, \dots, y_k \in \R,$ there exists a unique point $w \in
\R^d,$ such that the first $k$ coordinates of $w$ are $y_1, \dots,
y_k$ and $w$ is affinely dependent with the points in $U.$
\end{ilist}

\end{lem}

\begin{rem}\label{simplex}By \rm{(iv)}, the problem of counting
number of lattice points in a polytope in general position is
reduced to that of counting lattice points in a simplex in general
position. In particular, together with the fact that
$\vol(\bigsqcup_{i=1}^k P_i) = \sum_{i=1}^k \vol(P_i),$ to prove
Theorem \ref{main2} it is sufficient to prove the case when $P$ is a
lattice-face simplex.
\end{rem}

Therefore, we will only construct our decomposition in the case of
simplices in general position. However, before the construction, we
need one more proposition about the nonnegative part of a polytope
in general position.

\begin{prop}\label{psum}
Let $P$ be a $d$-polytope in general position with facets $F_1, F_2
\dots F_k.$ Let $H$ be the hyperplane determined by $F_k.$ For $i: 1
\le i \le k,$ let $F_i' = \pi^{-1}(\pi(F_i)) \cap H$ and $Q_i =
\conv(F_i \cup F_i').$

Then
\begin{equation}\label{esum}
\Omega(P) = - \sgn(F_k) \bigoplus_{i=1}^{k-1} \sgn(F_i)
\rho^+(\Omega(\pi(F_i)), Q_i).
\end{equation}

\end{prop}

\begin{proof}
We are going to just prove the case when $F_k$ is a negative facet;
for the other case we can prove it analogously. Suppose $F_1, F_2,
\dots, F_\ell$ are positive facets and $F_{\ell+1}, \dots, F_k$ are
negative facets.

A special case of Lemma \ref{lemone}/(iii) is when $R =
\Omega(\pi(P)),$ so we have $$\Omega(P) = \rho^+(\Omega(\pi(P)), P)
= \bigoplus_{y \in \Omega(\pi(P))} \rho^+(y, P).$$

Now for any points $a$ and $b,$ we use $(a, b]$ to denote the
half-open line segment between $a$ and $b$. Then, $\rho^+(y,P) =
(n(y,P), p(y,P)] = ( \rho(y,H), p(y, P)] \ominus (\rho(y, H),
n(y,P)] .$ Therefore,
\begin{align*}
\Omega(P) &=  \bigoplus_{y \in \Omega(\pi(P))} ( (\rho(y, H), p(y, P)] \ominus (\rho(y, H), n(y,P)] ) \\
&=\left(\bigoplus_{y \in \Omega(\pi(P))} (\rho(y, H), p(y,
P)]\right) \bigoplus \left(\bigoplus_{y \in \Omega(\pi(P))}
(-1)\cdot (\rho(y, H), n(y,P)] \right).
\end{align*}

By Lemma \ref{pgen}/(v), we have $\Omega(\pi(P)) =
\bigoplus_{i=1}^{\ell} \Omega(\pi(F_i)).$ Therefore,
\begin{eqnarray*}
\bigoplus_{y \in \Omega(\pi(P))} (\rho(y, H), p(y, P)] &=& \bigoplus_{i=1}^{\ell} \bigoplus_{y \in \Omega(\pi(F_i))} (\rho(y, H), p(y, P)] \\
&=& \bigoplus_{i=1}^{\ell} \bigoplus_{y \in \Omega(\pi(F_i))} (\rho(y, F_i'), \rho(y, F_i)] \\
&=& \bigoplus_{i=1}^{\ell} \rho^+(\Omega(\pi(F_i)), Q_i).
\end{eqnarray*}
Similarly, we will have
$$\bigoplus_{y \in \Omega(\pi(P))} (-1)\cdot (\rho(y, H), n(y,P)] = \bigoplus_{i=\ell+1}^k (-1) \rho^+(\Omega(\pi(F_i)), Q_i).$$
Note that $\rho^+(\Omega(\pi(F_k)), Q_k)$ is the empty set. Thus,
putting everything together, we get (\ref{esum}).

\end{proof}

Now, we can use this proposition to inductively construct a
decomposition of the nonnegative part $\Omega(P)$ of a $d$-simplex
$P$ in general position into $d!$ signed sets.

\subsection*{Decomposition of $\Omega(P)$:}
\begin{itemize}
\item If $d = 1,$ we do nothing: $\Omega(P) = \Omega(P).$
\item If $d \ge 2,$ then by applying Proposition \ref{psum} to $P$ and letting
$k = {d+1},$ we have
\begin{equation}\label{sum}
\Omega(P) = - \sgn(F_{d+1}) \bigoplus_{i=1}^d \sgn(F_i)
\rho^+(\Omega(\pi(F_i)), Q_i ).\end{equation}

However, by Lemma \ref{pgen}/{\rm (iii)}, each $\pi(F_i)$ is a
$(d-1)$-simplex in general position. By the induction hypothesis,
$\Omega(\pi(F_i)) = \bigoplus_{j=1}^{(d-1)!} S_{i,j},$ where
$S_{i,j}$'s are signed sets. $$\rho^+(\Omega(\pi(F_i)), Q_i ) =
\rho^+(\bigoplus_{j=1}^{(d-1)!} S_{i,j}, Q_i )=
\bigoplus_{j=1}^{(d-1)!} \rho^+( S_{i,j}, Q_i ).$$ Since each
$\rho^+( S_{i,j}, Q_i )$ is a signed set, we have decomposed
$\Omega(P)$ into $d!$ signed sets.
\end{itemize}

Now we know that we can decompose $\Omega(P)$ into $d!$ signed sets.
But we still need to figure out what these sets are and which signs
they have. In the next subsection, we are going to discuss the sign
of a facet of a $d$-simplex, which is going to help us determine the
signs in our decomposition.
\subsection{The sign of a facet of a $d$-simplex}\label{ssign}

 From now on, we will always
use the following setup for a $d$-simplex unless otherwise stated:

Suppose $P$ is a $d$-simplex in general position with vertex set $V
= \{v_1, v_2, \dots, v_{d+1}\},$ where the coordinates of $v_i$ are
$ \x_i = (x_{i,1},
x_{i,2}, \dots, x_{i,d}).$ 

For any $i,$ we denote by $F_i$ the facet determined by vertices in
$V \setminus \{v_i\}$ and $H_i$ the hyperplane determined by $F_i.$

For any $\sigma \in \S_d$ and $k: 1 \le k \le d,$ we define matrices
$X_V(\sigma, k)$ and $Y_V(\sigma, k)$ to be the matrices

$$X_V(\sigma, k) = \left(\begin{array}{ccccc}
1  & x_{\sigma(1), 1} & x_{\sigma(1),2} & \cdots & x_{\sigma(1), k} \\
1  & x_{\sigma(2), 1} & x_{\sigma(2),2} & \cdots & x_{\sigma(2), k} \\
\vdots  & \vdots & \vdots & \ddots & \vdots \\
1  & x_{\sigma(k), 1} & x_{\sigma(k),2} & \cdots & x_{\sigma(k), k} \\
1  & x_{d+1, 1} & x_{d+1,2} & \cdots & x_{d+1, k}
\end{array}\right)_{(k+1) \times (k+1)},$$

$$Y_V(\sigma, k) = \left(\begin{array}{ccccc}
1  & x_{\sigma(1), 1} & x_{\sigma(1),2} & \cdots & x_{\sigma(1), k-1} \\
1  & x_{\sigma(2), 1} & x_{\sigma(2),2} & \cdots & x_{\sigma(2), k-1} \\
\vdots  & \vdots & \vdots & \ddots & \vdots \\
1  & x_{\sigma(k), 1} & x_{\sigma(k),2} & \cdots & x_{\sigma(k),
k-1}
\end{array}\right)_{k \times k}.$$

We also define $z(\sigma, k)$ to be
$$z_V(\sigma,k) = \det(X_V(\sigma, k))/\det(Y_V(\sigma, k)),$$ where $\det(M)$ is the
determinant of a matrix $M.$

We often omit the subscript $V$ for $X_V(\sigma, k),$ $Y_V(\sigma,
k)$ and $z_V(\sigma,k)$ if there is no confusion.

Now we can determine the sign of a facet $F_i$ of $P$ by looking at
the determinants of these matrices, denoting by $\sgn(x)$ the usual
definition of sign of a real number $x.$

\begin{lem}\label{facetsign} We have
\begin{ilist}
\itm $\forall i: 1 \le i \le d$ and $\forall \sigma \in \S_d$ with
$\sigma(d) = i,$
\begin{equation}
\sgn(F_i) =  \sgn(\det(X(\sigma, d))/\det(X(\sigma, d-1))).
\end{equation}
\itm When $i = d+1$ and for $\forall \sigma \in \S_d,$
\begin{equation}
\sgn(F_{d+1}) = - \sgn(\det(X(\sigma, d))/\det(Y(\sigma, d)))= -
\sgn(z(\sigma, d)).
\end{equation}
\end{ilist}

\end{lem}

\begin{proof}
For any $i: 1 \le i \le d+1,$ let $v_i' = \rho(\pi(v_i), H_i),$ i.e.
$v_i'$ is the unique point of the hyperplane spanned by $F_i$ which
has the same coordinates as $v_i$ except for the last one. Suppose
the coordinates of $v_i'$ are $(x_{i,1}, \dots, x_{i, d-1},
x_{i,d}').$ Then $F_i$ is a positive facet if and only if $x_{i,d} <
x_{i,d}'.$ Therefore,
$$\sgn(F_i) = - \sgn(x_{i,d} - x_{i,d}').$$

$\forall i: 1 \le i \le d$ and $\forall \sigma \in \S_d$ with
$\sigma(d) = i,$ because $v_i'$ is in the hyperplane determined by
$F_i,$ we have that
$$\det\left(\left(\begin{array}{ccccc}
1  & x_{\sigma(1), 1} & \cdots & x_{\sigma(1),d-1} & x_{\sigma(1), d} \\
\vdots  & \vdots & \ddots & \vdots & \vdots \\
1  & x_{\sigma(d-1), 1} & \cdots & x_{\sigma(d-1),d-1} & x_{\sigma(d-1), d} \\
1  & x_{\sigma(d), 1} & \cdots & x_{\sigma(d),d-1} & x_{\sigma(d), d}' \\
1  & x_{d+1, 1} & \cdots & x_{d+1,d-1} & x_{d+1, d}
\end{array}\right)\right) = 0.$$
Therefore, $$ \det(X(\sigma, d)) = (-1)^{2d+1} (x_{i, d} - x_{i,d}')
\det(X(\sigma, d-1)).$$

Thus, $$\sgn(\det(X(\sigma, d))/\det(X(\sigma, d-1))) = -
\sgn(x_{i,d} - x_{i,d}') = \sgn(F_i).$$

We can similarly prove the formula for $i = d+1.$
\end{proof}

\subsection{Decomposition formulas}

The following theorem describes the signed sets in our
decomposition.

\begin{thm}\label{decomp} Let $P$ be a $d$-simplex in general position with vertex set
$V = \{v_1, v_2, \dots, v_{d+1}\},$ where the coordinates of $v_i$
are $ \x_i = (x_{i,1}, x_{i,2}, \dots, x_{i,d}).$ For any $\sigma
\in \S_d,$ and $k: 0 \le k \le d-1,$ let $v_{\sigma, k}$ be the
point with first $k$ coordinates the same as $v_{d+1}$ and affinely
dependent with $v_{\sigma(1)}, v_{\sigma(2)}, \dots, v_{\sigma(k)},
v_{\sigma(k+1)}.$ (By Lemma \ref{pgen}/{\rm (vii)}, we know that
there exists one and only one such point.) We also let $v_{\sigma,
d} = v_{d+1}.$ Then
\begin{equation}\label{fdecomp} \Omega(P) = \bigoplus_{\sigma
\in \S_d} \sgn(\sigma, P) S_{\sigma},
\end{equation}
 where
\begin{equation}\label{sign} \sgn(\sigma, P) = \sgn(\det(X(\sigma,
d))) \sgn\left(\prod_{i=1}^d z(\sigma, i)\right),
\end{equation}
and \begin{equation}\label{Ssigma} S_{\sigma}= \{ \s \in \R^d \ | \
\pi^{d-k}(\s) \in \Omega(\pi^{d-k}(\conv(\{v_{\sigma,0}, \dots,
v_{\sigma,k}\}))) \forall 1 \le k \le d\} \end{equation}
 is a set of
weight $1$, i.e. a regular set.

Hence,$$\L(\Omega(P)) = \bigoplus_{\sigma \in \S_d} \sgn(\sigma, P)
\L(S_{\sigma}).$$

\end{thm}
\begin{proof}
We prove it by induction on $d.$

When $d=1,$ the only permutation $\sigma \in \S_1$ is the identity
permutation $\1$. One can check that $\sgn(\1, P) = 1$ and
$S_{\1}=\Omega(\conv(v_1, v_2)).$ Thus (\ref{fdecomp}) holds.

Assuming (\ref{fdecomp}) holds for $d = d_0 \ge 1,$ we consider for
$d = d_0 + 1.$ For any $i: 1 \le i \le d,$ $\pi(F_i)$ is a
$(d-1)$-simplex in general position with vertex set $W = \{w_1,
\dots, w_d \},$ where $w_j = \begin{cases} \pi(v_j), & j < i, \\
\pi(v_{j+1}), & j \ge i. \end{cases}$ Therefore, by the induction
hypothesis,
\begin{equation}\Omega(\pi(F_i)) = \bigoplus_{\varsigma
\in \S_{d-1}} \sgn(\varsigma, \pi(F_i)) S_{\varsigma}',
\end{equation}
 where
\begin{equation*} \sgn(\varsigma, \pi(F_i)) = \sgn(\det(X_W(\varsigma,
d-1))) \prod_{i=1}^{d-1} \sgn(z_W(\varsigma, i)),
\end{equation*}
\begin{equation*} S_{\varsigma}'= \{ \s \in \R^{d-1} \ | \
\pi^{d-1-k}(\s) \in \Omega(\pi^{d-1-k}(\conv(\{w_{\varsigma,0},
\dots, w_{\varsigma,k}\}))) \forall 1 \le k \le d-1\}.
\end{equation*}
For any $\varsigma \in \S_{d-1},$ if we let $\sigma \in \S_d$ with
$\sigma(j) = \begin{cases}i, & j = d, \\ \varsigma(j), & \varsigma(j) < i, \\
\varsigma(j)+1, & \varsigma(j) \ge i, \end{cases}$ then this gives a
bijection between $\varsigma \in \S_{d-1}$ and $\sigma \in \S_d$
with $\sigma(d) = i.$ In particular, for any $j: 1 \le j \le d-1,$
$w_{\varsigma(j)} = \pi(v_{\sigma(j)}).$ Hence,
\begin{equation*} \sgn(\varsigma, \pi(F_i)) = \sgn(\det(X(\sigma,
d-1))) \prod_{i=1}^{d-1} \sgn(z(\sigma, i)).
\end{equation*}
Note that $w_{\varsigma, d-1} = w_d = \pi(v_{d+1}) = \pi(v_{\sigma,
d-1}),$ so
\begin{equation*} S_{\varsigma}'= \{ \s \in \R^{d-1} \ | \
\pi^{d-1-k}(\s) \in \Omega(\pi^{d-k}(\conv(\{v_{\sigma,0}, \dots,
v_{\sigma,k}\}))) \forall 1 \le k \le d-1\}.
\end{equation*}

One can check that $F_i' = \pi^{-1}(\pi(F_i) \cap H_{d+1} = \conv(\{
v_1, \dots, v_{i-1}, v_{i+1}, \dots, v_d, v_{\sigma, d-1}\})$ and
$Q_i = \conv(F_i \cup F_i') = \conv(V \cup \{v_{\sigma,d-1} \}
\setminus \{ v_i\} ).$ Hence,
$$\rho^+(S_{\varsigma}' , Q_i) = \{ \s \in \R^d \ | \
\pi^{d-k}(\s) \in \Omega(\pi^{d-k}(\conv(\{v_{\sigma,0}, \dots,
v_{\sigma,k}\}))) \forall 1 \le k \le d\}.$$

By letting $S_{\sigma} = \rho^+(S_{\varsigma}' , Q_i)$ and
$\sgn(\sigma, P) = - \sgn(F_{d+1}) \sgn(F_i) \sgn(\varsigma,
\pi(F_i))$ and using Lemma \ref{facetsign}, we get $$ -
\sgn(F_{d+1}) \sgn(F_i) \rho^+(\pi(F_i), Q_i) = \bigoplus_{\sigma
\in \S_d, \sigma(d) = i} \sgn(\sigma, P) S_{\sigma}.$$ Thus,
together with (\ref{sum}), summing over all $i: 1 \le i \le d$ gives
(\ref{fdecomp}).

\end{proof}

\begin{cor} If $P$ is a $d$-simplex in general position, then
\begin{equation}\label{lpf1}
|\L(\Omega(P))| = \sum_{\sigma \in \S_d} \sgn(\sigma, P)
|\L(S_{\sigma})|.
\end{equation}
\end{cor}

Therefore, if we can calculate the number of lattice points in
$S_{\sigma}$'s, then we can calculate the number of lattice points
in the nonnegative part of a $d$-simplex in general position.
Although it's not so easy to find $|\L(S_{\sigma}) |$'s for an
arbitrary polytope, we can do it for any lattice-face $d$-polytope.

\section{Lattice enumeration in $S_{\sigma}$ and Bernoulli polynomials}
In this section, we will count the number of lattice points in
$S_{\sigma}$'s when $P$ is a lattice-face $d$-simplex. This
calculation involves Bernoulli polynomials.

\subsection{Counting lattice points in $S_{\sigma}$}
We say a map from $\R^d \to \R^d$ is {\it lattice preserving} if it
is invertible and it maps lattice points to lattice points. Clearly,
given a lattice preserving map $f,$  for any set $S \in \R^d$ we
have that $|\L(S)| = |\L(f(S))|.$

Let $P$ be a lattice face $d$-simplex with vertex set $V =\{v_1,
\dots, v_{d+1} \},$ where we use the same setup as before for
$d$-simplices.

Given any $\sigma \in \S_d,$ recall that $S_{\sigma}$ is defined as
in (\ref{Ssigma}). To count the number of lattice points in
$S_{\sigma},$ we want to find a lattice preserving affine
transformation which simplifies the form of $S_{\sigma}.$

Before trying to find such a transformation, we will define more
notation.

For any $\sigma \in \S_d,$ $k: 1 \le k \le d$ and $\x = (x_1, x_2,
\dots, x_d) \in \R^d,$ we define matrix $\tilde{X}(\sigma, k; \x)$
as

$$\tilde{X}(\sigma, k; \x) = \left(\begin{array}{ccccc}
1  & x_{\sigma(1), 1} & x_{\sigma(1),2} & \cdots & x_{\sigma(1), k} \\
1  & x_{\sigma(2), 1} & x_{\sigma(2),2} & \cdots & x_{\sigma(2), k} \\
\vdots  & \vdots & \vdots & \ddots & \vdots \\
1  & x_{\sigma(k), 1} & x_{\sigma(k),2} & \cdots & x_{\sigma(k), k} \\
1  & x_{1} & x_{2} & \cdots & x_{k}
\end{array}\right)_{(k+1) \times (k+1)},$$
and for $j: 0 \le j \le k,$ let $\m(\sigma, k; j)$ be the minor of
the matrix $\tilde{X}(\sigma, k; \x)$ obtained by omitting the last
row and the $(j+1)$th column. Then
\begin{equation}\label{detm1}
\det(\tilde{X}(\sigma, k; \x)) = (-1)^k \left(\m(\sigma, k; 0) +
\sum_{j=1}^k (-1)^j  \m(\sigma, k; j) x_j\right).
\end{equation}

Note that $\m(\sigma, k; k) = \det(Y(\sigma, k)).$ Therefore,
\begin{equation}\label{detm2}
\frac{\det(\tilde{X}(\sigma, k; \x))}{\det(Y(\sigma, k))} = (-1)^k
\frac{\m(\sigma, k; 0)}{\det(Y(\sigma,k))} + \sum_{j=1}^{k-1}
(-1)^{k+j} \frac{\m(\sigma, k; j)}{\det(Y(\sigma,k))}x_j + x_k.
\end{equation}

We will construct our transformation based on (\ref{detm2}). Before
that, we give the following lemma which discusses the coefficients
in the right hand side of (\ref{detm2}).

\begin{lem}\label{int1}
Suppose $P$ is a lattice-face $d$-simplex. $\forall \sigma \in \S_d,
\forall k: 1 \le k \le d,$ and $\forall j: 0 \le j \le k-1,$ we have
that
$$\frac{\m(\sigma, k; j)}{\det(Y(\sigma,k))} \in \Z.$$
\end{lem}

\begin{proof}
By the definition of lattice-face polytope and Lemma \ref{plf}/{\rm
(i)}, one can see that $\pi^{d-k}(\conv(v_{\sigma(1)}, \dots,
v_{\sigma(k)}, v_{d+1})) = \conv(\pi^{d-k}(v_{\sigma(1)}), \dots,
\pi^{d-k}(v_{\sigma(k)}), \pi^{d-k}(v_{d+1}))$ is a lattice-face
$k$-polytope. Choose $U = \{\pi^{d-k}(v_{\sigma(1)}), \dots,
\pi^{d-k}(v_{\sigma(k)}) \},$ then $\pi(\L(H_U)) = \Z^{k-1},$ where
$H_U$ is the affine space spanned by $U.$ However,
$$H_U = \{ \x = (x_1, \dots, x_k) \in \R^k \ | \ \det(\tilde{X}(\sigma, k; \x))
= 0 \}.$$ Therefore, we must have that
$$\det(\tilde{X}(\sigma, k; \x)) = 0, \quad x_1, \dots, x_{k-1} \in \Z
\Rightarrow x_k \in Z.$$

Let $x_1 = \cdots = x_{k-1} = 0,$ then $\det(\tilde{X}(\sigma, k;
\x)) = 0$ implies that $$(-1)^k \frac{\m(\sigma, k;
0)}{\det(Y(\sigma,k))} + x_k = 0 \Rightarrow \frac{\m(\sigma, k;
0)}{\det(Y(\sigma,k))} = (-1)^{k+1} x_k \in Z.$$

For any $j: 1 \le j \le k-1,$ let $x_i = x_{\sigma(j),i} +
\delta_{i,j}$ for $1 \le i \le k-1,$ where $\delta_{i,j}$ is the
Kronecker delta function. Then, $\det(\tilde{X}(\sigma,k; x)) = 0$
implies that
\begin{eqnarray*}
0 &=& (-1)^k \frac{\m(\sigma, k; 0)}{\det(Y(\sigma,k))} +
\sum_{i=1}^{k-1} (-1)^{k+i} \frac{\m(\sigma, k;
i)}{\det(Y(\sigma,k))}x_{\sigma(j),i} + x_k + (-1)^{k+j}
\frac{\m(\sigma, k; j)}{\det(Y(\sigma,k))} \\
&=&x_k - x_{\sigma(j),k} + (-1)^{k+j} \frac{\m(\sigma, k;
j)}{\det(Y(\sigma,k))},
\end{eqnarray*}
where the second equality follows from the fact that
$(x_{\sigma(j),1}, \dots, x_{\sigma(j), k})$ is in $H_U.$ Thus,
$$\frac{\m(\sigma, k;
j)}{\det(Y(\sigma,k))} = (-1)^{k+j+1}(x_k - x_{\sigma(j),k}) \in
\Z.$$

\end{proof}

Given this lemma, we have the following proposition.
\begin{prop}
There exist a lattice-preserving affine transformation $T_{\sigma}$
which maps $x = (x_1, x_2, \dots, x_d) \in \R^d$ to
$$(\frac{\det(\tilde{X}(\sigma, 1; \x))}{\det(Y(\sigma,1))},
\frac{\det(\tilde{X}(\sigma, 2; \x))}{\det(Y(\sigma,2))}, \dots,
\frac{\det(\tilde{X}(\sigma, d; \x))}{\det(Y(\sigma, d))}).$$
\end{prop}

\begin{proof}
Let $\alpha_{\sigma} = (- \frac{\m(\sigma, 1;
0)}{\det(Y(\sigma,1))}, \frac{\m(\sigma, 2; 0)}{\det(Y(\sigma,2))},
\dots, (-1)^d \frac{\m(\sigma, d; 0)}{\det(Y(\sigma,d))})$ and
$M_{\sigma} = (m_{\sigma, j, k})_{d \times d},$ where

$$m_{\sigma, j, k} = \begin{cases}1, &\mbox{if $j = k$}  \\ 0, &\mbox{if $j > k$} \\
(-1)^{k+j}\frac{\m(\sigma, k; j)}{\det(Y(\sigma,k))}, &\mbox{if $j <
k$}\end{cases}.$$

We define $T_{\sigma}: \R^d \to \R^d$ by mapping $\x$ to
$\alpha_{\sigma} + \x M_{\sigma}.$ By (\ref{detm2}),
$$\alpha_{\sigma} + \x M_{\sigma} = (\frac{\det(\tilde{X}(\sigma, 1;
\x))}{\det(Y(\sigma,1))}, \frac{\det(\tilde{X}(\sigma, 2;
\x))}{\det(Y(\sigma,2))}, \dots, \frac{\det(\tilde{X}(\sigma, d;
\x))}{\det(Y(\sigma, d))}).$$ Also, because all of the entries in
$M_{\sigma}$ and $\alpha_{\sigma}$ are integers and the determinant
of $M_{\sigma}$ is $1,$ $T_{\sigma}$ is lattice preserving.
\end{proof}

\begin{cor}\label{coord} Give $P$ a lattice-face polytope with vertex set $V =
\{v_1, v_2, \dots, v_{d+1} \},$ we have that
\begin{ilist}
\itm $\forall i: 1 \le i \le d,$ the last $d+1-i$ coordinates of
$T_{\sigma}(v_{\sigma(i)})$ are all zero.

\itm $$T_{\sigma}(v_{d+1}) = (z(\sigma,1), z(\sigma,2), \dots,
z(\sigma,d)).$$

\itm Recall that for $k: 0 \le k \le d-1$, $v_{\sigma, k}$ is the
unique point with first $k$ coordinates the same as $v_{d+1}$ and
affinely dependent with $v_{\sigma(1)}, v_{\sigma(2)}, \dots,
v_{\sigma(k)}, v_{\sigma(k+1)}.$ Then the first $k$ coordinates of
$T_{\sigma}(v_{\sigma, k})$ are the same
 as $T_{\sigma}(v_{d+1})$ and the rest of the coordinates are zero.
 In other words, $$T_{\sigma}(v_{\sigma,k}) = (z(\sigma,1), \dots,
 z(\sigma,k), 0, \dots, 0).$$
\end{ilist}
\end{cor}

\begin{proof}
{\rm (i)} This follows from that fact that $\det(\tilde{X}(\sigma,
k; \x_{\sigma(i)})) = 0$ if $1 \le i \le k \le d.$

{\rm (ii)} This follows from the fact that $\tilde{X}(\sigma, k;
\x_{d+1}) = {X}(\sigma, k)$ and $z(\sigma,k) = \det(X(\sigma,
k))/\det(Y(\sigma, k)).$

{\rm (iii)} Because for any $\x \in \R^d,$ the $k$th coordinate of
$T_{\sigma}(\x)$ only depends on the first $k$ coordinates of $\x,$
$T_{\sigma}(v_{\sigma, k})$ has the same first $k$ coordinates as
$T_{\sigma}(v_{d+1}).$ $T_{\sigma}$ is an affine transformation. So
$T_{\sigma}(v_{\sigma, k})$ is affinely dependent with
$T_{\sigma}(v_{\sigma(1)}),$ $T_{\sigma}(v_{\sigma(2)}),$ $\dots,$
$T_{\sigma}(v_{\sigma(k)}),$ $T_{\sigma}(v_{\sigma(k+1)}),$ the last
$d-k$ coordinates of which are all zero. Therefore the last $d-k$
coordinates of $T_{\sigma}(v_{\sigma, k})$ are all zero as well.
\end{proof}

Recalling that $v_{\sigma, d} = v_{d+1},$ we are able to describe
$T_{\sigma}(S_{\sigma})$ now.
\begin{prop}
Let $\nHat{S}_{\sigma} = T_{\sigma}(S_{\sigma}).$ Then
\begin{equation}\label{Ssigmahat}
\s = (s_1, s_2, \dots, s_d) \in \nHat{S}_{\sigma}  \Leftrightarrow
\forall 1 \le k \le d, s_k \in \Omega(\conv(0,
\frac{z(\sigma,k)}{z(\sigma,k-1)}s_{k-1} )),
\end{equation}
where by convention we let $z(\sigma, 0) = 1$ and $s_0 = 1.$

\end{prop}
\begin{proof} $T_{\sigma}$ is an affine transformation whose
corresponding matrix $M_{\sigma}$ is upper triangular. So
$T_{\sigma}$ commutes with $\Omega,$ $\pi$ and $\conv.$ Therefore,
$$\nHat{S}_{\sigma}= \{ \s \in \R^d \ | \
\pi^{d-k}(\s) \in \Omega(\pi^{d-k}(\conv(\{\nHat{v}_{\sigma,0},
\dots, \nHat{v}_{\sigma,k}\}))) \forall 1 \le k \le d\},$$
 where
$\nHat{v}_{\sigma,i} = T_{\sigma}(v_{\sigma,i}) = (z(\sigma, 1),
\dots, z(\sigma,i), 0, \dots, 0),$ for $0 \le i \le d.$

(\ref{Ssigmahat}) follows.
\end{proof}

Because $T_{\sigma}$ is a lattice preserving map, $|\L(S_{\sigma})|
= |\L(\nHat{S}_{\sigma})|.$ Hence, our problem becomes to find the
number of lattice points in $\nHat{S}_{\sigma}$. However,
$\nHat{S}_{\sigma}$ is much nicer than ${S}_{\sigma}.$ Actually, we
can give a formula to calculate all of the sets having the same
shape as $\nHat{S}_{\sigma}.$

\begin{lem}\label{lfS} Given real nonzero numbers $b_0 = 1, b_1, b_2, \dots, b_d,$ let $a_k' = b_k/b_{k-1}$ and $a_k = b_k/|b_{k+1}|, \forall k : 1 \le k \le d.$ Let $S$ be the
set defined by the following:
$$\s = (s_1, s_2, \dots, s_d) \in S \Leftrightarrow \forall 1 \le k
\le d, s_k \in \Omega(\conv(0, a_k' s_{k-1} )),$$ where $s_0$ is set
to $1.$ Then \begin{equation}\label{lfS1} |\L(S)| = \sum_{s_1 \in
\L(\Omega(\conv(0, a_1' )))} \sum_{s_2 \in \L(\Omega(\conv(0, a_2'
s_{1} )))} \cdots \sum_{s_d \in \L(\Omega(\conv(0, a_d' s_{d-1} )))}
1.
\end{equation}
 In particular, if $b_d > 0,$ then
 \begin{equation}\label{lfS2}
|\L(S)| = \sum_{s_1 = 1}^{\nBar{\lfloor a_1 \rfloor}} \sum_{s_2
=1}^{\nBar{\lfloor a_2 s_1 \rfloor}} \cdots \sum_{s_d
=1}^{\nBar{\lfloor a_d s_{d-1} \rfloor}} 1, \end{equation} where for
any real number $x,$ $\lfloor x \rfloor$ is the largest integer no
greater than $x$ and $\nBar{x}$ is defined as $$\nBar{x} =
\begin{cases}x, & \mbox{if $x \ge 0$} \\ -x-1, &\mbox{if
$x < 0$}
\end{cases}.$$

\end{lem}

Note that $\nBar{\lfloor x \rfloor} \in \Z_{\ge 0}$, and if any of
the sums in (\ref{lfS2}) have upper bound equal to $0,$ we consider
the sum to be $0.$

\begin{proof}
(\ref{lfS1}) is straightforward. (\ref{lfS2}) follows from the facts
that for any real numbers $x,$
$$\L(\Omega(\conv(0,x))) =
\begin{cases}\{z \in \Z \ | \ 1 \le z \le \nBar{\lfloor x \rfloor}
\} &\mbox{if $x \ge 0$}\\
\{z \in \Z \ | \ -\nBar{\lfloor x \rfloor} \le z \le 0 \} &\mbox{if
$x < 0$}
\end{cases},$$
the sign of $s_i$ is the same as the sign of $b_i,$ and, because
$b_d
> 0,$ all the $s_i$'s are non-zero.
\end{proof}

We want to give a formula for the number of lattice points in
$\nHat{S}_{\sigma}$ in the form of (\ref{lfS2}). We first need the
condition ``$b_d > 0$'', which in our case is that ``$z(\sigma, d)
> 0$''. However, for any $d$-simplex $P$ in general position, we can always
find a way to order its vertices into $V = \{v_1, v_2, \dots,
v_{d+1} \},$ so that the corresponding $\det(X(\1, d))$ and
$\det(Y(\1, d))$ are positive, where $\1$ stands for the identity
permutation in $\S_d.$ Note $z(\sigma, d)$ is independent of
$\sigma.$ So it is positive.

Moreover, for lattice-polytopes, we have another good property of
the $z(\sigma, k)$'s which allows us to remove the $\lfloor \
\rfloor$ operation in (\ref{lfS2}).

\begin{lem}\label{zsigma} If $P$ is a lattice-polytope $d$-simplex, then
$$z(\sigma,k)/z(\sigma, k-1) \in \Z,$$
where by convention $z(\sigma,0)$ is set to $1.$
\end{lem}

\begin{proof}
Let $P' = T_{\sigma}(P)$ with vertex set $V' = \{v_1', \dots,
v_{d+1}'\},$ where $v_i' = T_{\sigma}(v_i)$ with coordinates $\x_i'
= (x_{i,1}', \dots, x_{i,d}').$ Because $T_\sigma$ is a upper
triangular lattice preserving map, $P'$ is a lattice-face
$d$-simplex as well. Similar to the proof of Lemma \ref{int1},
$\conv(\pi^{d-k}(v_{\sigma(1)}'), \dots, \pi^{d-k}(v_{\sigma(k)}'),
\pi^{d-k}(v_{d+1}'))$ is a lattice-face $k$-polytope. We choose $U =
\{\pi^{d-k}(v_{\sigma(1)}'), \dots, \pi^{d-k}(v_{\sigma(k-1)}'),
\pi^{d-k}(v_{d+1}')\},$ then $\pi(\L(H_U)) = \Z^{k-1}.$ Note that by
Corollary \ref{coord}/{\rm (i)},{\rm (ii)}, we have that
\begin{alist}
\itm the last $2$ coordinates of $\pi^{d-k}(v_{\sigma(j)}')$ are
both zero, for any $j: 1 \le j \le k-1.$

\itm $\pi^{d-k}(v_{d+1}') = (z(\sigma, 1), \dots, z(\sigma, k-1),
z(\sigma, k)).$
\end{alist}
Hence, $(x_1, \dots, x_k) \in H_U$ if and only if
$\det\left(\left(\begin{array}{cc} z(\sigma, k-1) & z(\sigma, k) \\
x_{k-1} & x_k \end{array}\right)\right) = 0,$ where we set $x_0 =
1.$

We have that for any $(x_1, \dots, x_k) \in H_U,$ if $x_1, \dots,
x_{k-1} \in \Z,$ then $x_k \in \Z.$ Thus, by setting $x_{k-1} = 1,$
we get $z(\sigma,k)/z(\sigma, k-1) = x_k \in \Z.$

\end{proof}

Therefore, by Lemma \ref{lfS} and Lemma \ref{zsigma}, we have the
following result.

\begin{prop} Let $P$ be a lattice-face $d$-simplex with vertex set
$V,$ where the order of vertices makes both $\det(X(\1, d))$ and
$\det(Y(\1, d))$ positive. Define $$a(\sigma, k) =
\frac{z(\sigma,k)}{|z(\sigma,k-1)|}, \forall k: 1 \le k \le d.$$
Then
\begin{equation}\label{lSsigma1}
|\L(S_{\sigma})| = \sum_{s_1 = 1}^{\nBar{a(\sigma,1)}} \sum_{s_2
=1}^{\nBar{ a(\sigma,2) s_1 }} \cdots \sum_{s_d =1}^{\nBar{
a(\sigma,d) s_{d-1} }} 1.\end{equation}
\end{prop}

Because of (\ref{lSsigma1}), it's natural for us to define
\begin{equation}\label{dfd}
f_d(a_1, a_2, \dots, a_d) = \sum_{s_1 =1}^{a_1} \sum_{s_2 =1}^{a_2
s_1} \cdots \sum_{s_d =1}^{a_d s_{d-1}} 1,
\end{equation}
for any positive integers $a_1, a_2, \dots, a_d.$

Also, fixing $b_0 = 1,$ we define
\begin{equation}\label{dgd}
g_d(b_1, b_2, \dots, b_d) = f_d(b_1/b_0, b_2/b_1, \dots,
b_d/b_{d-1}), \end{equation} for any $(b_1, b_2, \dots, b_d) \in
(\Z_{>0})^d$ such that $b_i$ is a multiple of $b_{i-1} (\forall 1
\le i \le d).$


$f_d$ and $g_d$ are closely related to formula (\ref{lSsigma1}). In
next subsection, we will discuss Bernoulli polynomials and power
sums, which are connected to $f_d$ and $g_d,$ and then rewrite
(\ref{lSsigma1}) in terms of $g_d.$

\subsection{Power sums and Bernoulli polynomials}


The {\it $k$th Bernoulli polynomial}, $B_k(x)$, is defined as
\cite[p. 804]{AbrSte}
$$\frac{t e^{tx}}{e^t -1} = \sum_{k=0}^{\infty} B_k(x)
\frac{t^k}{k!},$$

The Bernoulli polynomials satisfy \cite{Lehmer}
\begin{equation}\label{bp1}
B_k(1-x) = (-1)^k B_k(x), \forall k \ge 0, \end{equation}
 as well as the relation \cite[p. 127]{WhiWat}
\begin{equation}\label{bp2}
B_k(x+1) - B_k(x) = k x^{k-1}, \forall k \ge 1. \end{equation} We
call $B_k = B_k(0)$ a {\it Bernoulli number}. It satisfies
\cite{Bnumber} that
\begin{equation}\label{bn}
B_k(0) = 0, \mbox{for any odd number $k \ge 3.$}\end{equation}

For $k \ge 0,$ let $$P_k(x) = \frac{B_{k+1}(x+1) - B_{k+1}}{k+1}.$$
Given any $n$ a positive integer, by (\ref{bp2}), we have that
$$P_k(n) = \sum_{i=0}^n i^k = \begin{cases}
\sum_{i=1}^n i^k & \mbox{if $k \ge 1$}\\
n + 1 & \mbox{if $k = 0$}
\end{cases}.$$  Therefore, we call $P_k(x)$ the {\it
$k$th power sum polynomial.} It's well known that for $k \ge 1,$
\begin{eqnarray}
& &\mbox{$P_k(x)$ is a polynomial in $x$ of degree $k+1,$} \label{pdeg}\\
& &\mbox{the constant term of $P_k(x)$ is $0,$ i.e., $x$ is a factor
of $P_k(x),$} \label{pconst}\\
& &\mbox{the leading coefficient of $P_k(x)$ is $\frac{1}{k+1}.$}
\label{plead}
\end{eqnarray}

\begin{lem} For any $k \ge 1,$
\begin{equation}\label{recp}
P_k(x) = (-1)^{k+1} P_k(-x-1).
\end{equation}
\end{lem}
\begin{proof}
It follows from (\ref{bp1}) and (\ref{bn}).
\end{proof}

\subsection*{Extension of the sum operation} Given $h = h(s)=\sum_{k \ge 0} h_k
s^k$ a polynomial in $s,$ the upper bound $u$ of a sum $\sum_{s=1}^u
h$ should be a positive integer in the usual definition. We extend
this definition to allow $u$ (as well as the $h_k$'s) to be in any
polynomial ring over $\R$ using the formula
\begin{equation}\label{extend}
\sum_{s=1}^u h = h_0 u + \sum_{k \ge 1} h_k P_k(u).
\end{equation}
One can check that this extension agrees with the case when $u$ is a
positive integer.

Since $f_d$ is defined by (\ref{dfd}), which recursively uses the
sum operation, we can use (\ref{extend}) to extend the domain of
$f_d$ from $\Z_{> 0}^d$ to $\Z^d$ or even $\R^d.$ Hence, the domain
of $g_d$ can be extended to $(\R \setminus \{0\})^d.$

\begin{lem}
$f_d(a_1, \dots, a_d)$ is a polynomial in $a_1$ of degree $d.$ And
$\prod_{i=1}^d a_i$ is a factor of it. In particular, $f_d$ can be
written as
\begin{equation}\label{fd}
f_d(a_1, \dots, a_d) = \sum_{k=1}^d f_{d,k}(a_2, \dots, a_d) a_1^k,
\end{equation}
where $f_{d,k}(a_2, \dots, a_d)$ is a polynomial in $a_2, \dots,
a_d$ with a factor $\prod_{i=2}^d a_i.$
\end{lem}
\begin{proof}
This can be proved by induction on $d$, using (\ref{pdeg}) and
(\ref{pconst}).
\end{proof}

\begin{lem}\label{rmbar1}
Given $(a_1, a_2, \dots, a_d) \in \R^d,$
$$f_d(a_1,a_2, \dots, a_d) = - \sum_{s_1=1}^{-a_1-1} \sum_{s_2=1}^{-a_2 s_1}
\sum_{s_3=1}^{a_3 s_2}\cdots \sum_{s_d =1}^{a_{d} s_{d-1}} 1.$$
\end{lem}

\begin{proof}
%
$$- \sum_{s_1=1}^{-a_1-1} \sum_{s_2=1}^{-a_2 s_1}
\sum_{s_3=1}^{a_3 s_2}\cdots \sum_{s_d =1}^{a_{d} s_{d-1}} 1 = -
\sum_{s_1=1}^{-a_1-1}  f_{d-1}(-a_2 s_1, a_3, \dots, a_d).$$ By
(\ref{fd}) and (\ref{recp}), we have
$$f_{d-1}(-a_2 s_1, a_3, \dots, a_d) = \sum_{k=1}^{d-1} f_{d-1,k}(a_3, \dots,
a_d) (-a_2 s_1)^k,$$ and
$$\sum_{s_1=1}^{-a_1-1} s_1^k = P_k(-a_1-1) =
(-1)^{k+1}P_k(a_1) = (-1)^{k+1} \sum_{s_1=1}^{a_1} s_1^k.$$

Therefore,

\begin{eqnarray*}
& &- \sum_{s_1=1}^{-a_1-1} \sum_{s_2=1}^{-a_2 s_1} \sum_{s_3=1}^{a_3
s_2}\cdots \sum_{s_d =1}^{a_{d} s_{d-1}} 1 = \sum_{s_1}^{a_1}
\sum_{k=1}^{d-1} f_{d-1,k}(a_3, \dots, a_d)
(a_2 s_1)^k \\
&=& \sum_{s_1}^{a_1} f_{d-1}(a_2 s_1, a_3, \dots, a_d) = f_d(a_1,
a_2, \dots, a_d).
\end{eqnarray*}

\end{proof}

\begin{prop}\label{rmbar2}
Given $b_0 = 1,$ $b_1, b_2, \dots, b_d \in (\R \setminus \{0\})$
with $b_d
> 0,$ let $a_k = b_k/|b_{k-1}|,$ then

\begin{equation}\label{gd}
g_d(b_1,b_2, \dots, b_d) = f_d(b_1, \frac{b_2}{b_1}, \dots,
\frac{b_d}{b_{d-1}}) = \sgn\left( \prod_{i=1}^d b_i \right)\sum_{s_1
= 1}^{\nBar{a_1}} \sum_{s_2 =1}^{\nBar{ a_2 s_1 }} \cdots \sum_{s_d
=1}^{\nBar{ a_d s_{d-1} }} 1,\end{equation} where we always treat
$s_i$ as positive when determining the meaning of $\nBar{ a_{i+1}
s_i }.$ That is, for $a_{i+1} > 0,$ we set $\nBar{ a_{i+1} s_i } =
a_{i+1} s_i,$ and for $a_{i+1} < 0,$ we set $\nBar{ a_{i+1} s_i } =
-a_{i+1} s_i-1.$ Note that this agrees with the original definition
when the $a_i$'s are all positive integers.
\end{prop}

\begin{proof}
We prove the proposition by induction on $d.$ When $d=1,$ it's
trivial.

Assume (\ref{gd}) holds for $d = d_0 \ge 1.$ $s_1$ is positive.
Thus, by the induction hypothesis,
$$g_d(\frac{b_2}{|b_1|}s_1, \frac{b_3}{|b_1|}s_1, \dots,
\frac{b_d}{|b_1|}s_1) = f_d(\frac{b_2}{|b_1|}s_1 ,\frac{b_3}{b_2},
\dots, \frac{b_d}{b_{d-1}}) = \sgn\left( \prod_{i=2}^d b_i \right)
\sum_{s_2 =1}^{\nBar{ a_2 s_1 }} \cdots \sum_{s_d =1}^{\nBar{ a_d
s_{d-1} }} 1 .$$ It's clear that (\ref{gd}) holds when $b_1 > 0.$ In
the case that $b_1 < 0,$ (\ref{gd}) follows from the above equation
and Lemma \ref{rmbar1}.
\end{proof}

\begin{prop} Let $P$ be a lattice-face $d$-simplex with vertex set
$V,$ where the order of vertices makes both $\det(X(\1, d))$ and
$\det(Y(\1, d))$ positive. Then
\begin{equation}\label{lSsigma2}
|\L(S_{\sigma})| = \sgn\left(\prod_{i=1}^d z(\sigma,i)\right)
g_d(z(\sigma,1), z(\sigma,2), \dots, z(\sigma,d)).\end{equation}
Therefore,
\begin{equation}\label{lpf2}
|\L(\Omega(P))| =  \sum_{\sigma \in \S_d}
\sgn(\sigma)g_d(z(\sigma,1), z(\sigma,2), \dots, z(\sigma,d)).
\end{equation}
\end{prop}

\begin{proof}
We can get (\ref{lSsigma2}) by comparing (\ref{lSsigma1}) and
(\ref{gd}). And (\ref{lpf2}) follows from (\ref{lpf1}),
(\ref{sign}), (\ref{lSsigma2}) and the fact that $\det(X(\sigma, d))
= \sgn(\sigma) \det(X(\1, d)).$
\end{proof}

\section{Proof of the Main Theorems}

We now have all the ingredients but one to prove our main theorems:
Theorem \ref{main1} and Theorem \ref{main2}. The missing one is
stated as the following proposition and will be proved in the next
section.

\begin{prop}\label{gdsigma}Let $V = \{v_1, v_2, \dots,
v_{d+1}\}$ be the vertex set of a $d$-simplex in general position,
where the coordinates of $v_i$ are $ \x_i = (x_{i,1}, x_{i,2},
\dots, x_{i,d}).$ Recall that $X(\sigma,k),$ $Y(\sigma,k)$ and
$z(\sigma,k)$ are defined in \textsection \ref{ssign} and $g_d$ is
defined in (\ref{dgd}). Then
\begin{equation}\label{gsigma}
\sum_{\sigma \in \S_d} \sgn(\sigma)g_d(z(\sigma,1), z(\sigma,2),
\dots, z(\sigma,d)) = \frac{1}{d!}\det(X(\1, d)),
\end{equation}
where $\1$ is the identity in $\S_d.$
\end{prop}

Given this proposition, we can prove the theorems.

\begin{proof}[Proof of Theorem \ref{main1} and Theorem \ref{main2}]
 As we mentioned in Remark \ref{simplex}, to prove Theorem
\ref{main2}, it is sufficient to prove the case when $P$ is a
lattice-face simplex.

When $P$ is a lattice-face $d$-simplex, we still assume that the
order of the vertices of $P$ makes both $\det(X(\1, d))$ and
$\det(Y(\1, d))$ positive. Thus, (\ref{lpf2}), (\ref{gsigma}) and
the fact that the volume of $P$ is $\frac{1}{d!}|\det(X(\1, d))|$
imply Theorem \ref{main2}, and Theorem \ref{main1} follows.
\end{proof}

Recall that we use $I(P)$ to denote the interior of  a $d$-polytope
$P.$ We denote by $\nHat{i}(P,m) = |I(m P) \cap \Z^d|$ the number of
lattice points in the interior of $m P.$

\begin{cor} For any lattice-face $d$-polytope $P,$ we have that
\begin{equation}\label{interior}
\nHat{i}(P, m) = \vol(m P) - \nHat{i}(\pi(P), m) = \sum_{k=0}^d
(-1)^{d-k} \vol_k(\pi^{(d-k)}(P)) m^k.
\end{equation}
Thus,
\begin{equation}\label{reciprocity}
i(P, -m) = (-1)^d \nHat{i}(P,m).
\end{equation}
\end{cor}

\begin{proof}
Since $P$ satisfies (\ref{one}), by Lemma \ref{lemone}/{\rm (vi)}
and Lemma \ref{plf}/{\rm (vi)}, $\pi$ induces a bijection between
$\L(PB(P) \cap NB(P))$ and $\L(\partial \pi(P)).$ Together with
Lemma \ref{plf}/{\rm (ii), (iii)}, this implies
$$\nHat{i}(P,m) = i(P,m) -  i(\pi(P),m) - \nHat{i}(\pi(P),m).$$
Therefore, (\ref{interior}) and (\ref{reciprocity}) follow from
Theorem \ref{main1}.
\end{proof}

Note that (\ref{reciprocity}) recovers the Ehrhart-Macdonald
reciprocity law \cite{Macdonald}, which states that for any integral
$d$-polytope, (\ref{reciprocity}) holds.

The proof of Proposition \ref{gdsigma} is self-contained and
different from the rest of this paper. We put it separately in the
next section.

\section{Proof of Proposition \ref{gdsigma}}
The purpose of this section is to prove Proposition \ref{gdsigma} by
showing both sides of (\ref{gdsigma}) are equal to
$$\frac{1}{d!} \sum_{\sigma \in \S_d} \sgn(\sigma)
\prod_{j=1}^d z(\sigma, j).$$ We always assume that $V = \{v_1, v_2,
\dots, v_{d+1}\}$ is the vertex set of a $d$-simplex in general
position, where the coordinates of $v_i$ are $ \x_i = (x_{i,1},
x_{i,2}, \dots, x_{i,d}).$ We first give some new notation and
definitions.

For all $k: 1 \le k \le d,$ let $\nHat{\x}_k = (\nHat{x}_{k,1},
\dots, \nHat{x}_{k,d}) = (x_{k,1}-x_{d+1,1}, \dots, x_{k,d} -
x_{d+1,d}) = \x_k - \x_{d+1}.$ Define
$$ \nHat{X}_V(\sigma,k) = \left(\begin{array}{ccccc}
\nHat{x}_{\sigma(1), 1} & \nHat{x}_{\sigma(1),2} & \cdots & \nHat{x}_{\sigma(1), k} \\
\nHat{x}_{\sigma(2), 1} & \nHat{x}_{\sigma(2),2} & \cdots & \nHat{x}_{\sigma(2), k} \\
 \vdots & \vdots & \ddots & \vdots \\
\nHat{x}_{\sigma(k), 1} & \nHat{x}_{\sigma(k),2} & \cdots &
\nHat{x}_{\sigma(k), k}
\end{array}\right),$$
$$\nHat{Y}_V(\sigma, k) = \left(\begin{array}{ccccc}
1  & \nHat{x}_{\sigma(1), 1}  & \cdots & \nHat{x}_{\sigma(1), k-1} \\
1  & \nHat{x}_{\sigma(2), 1} & \cdots & \nHat{x}_{\sigma(2), k-1} \\
\vdots  & \vdots  & \ddots & \vdots \\
1  & \nHat{x}_{\sigma(k), 1}  & \cdots & \nHat{x}_{\sigma(k), k-1}
\end{array}\right),$$ and
$$\nHat{z}_V(\sigma,k)= \det(\nHat{X}_V(\sigma,k))/\det(\nHat{Y}_V(\sigma,k)).$$
Then \begin{equation}\label{zhat} \nHat{z}_V(\sigma,k)= (-1)^k
z_V(\sigma,k).\end{equation}

Again, when there is no confusion, we omit the subscript $V$ from
$\nHat{X}_V(\sigma,k),$ $\nHat{Y}_V(\sigma,k)$ and
$\nHat{z}_V(\sigma,k).$

We define certain subsets of the symmetric group $\S_d,$ which we
will use in our later proofs. We denote by $\S_T$ the set of
permutations on some set $T$ and use one-line notation for all
permutations.
\begin{defn}
\begin{alist}
\itm Let $(\Lambda, \Gamma, \Delta)$ be a partition of $[d]$ with
the sizes of $\Lambda$ and $\Gamma$ to be $\ell$ and $i,$
respectively. For any $\lambda \in \S_{\Lambda},$ $\gamma \in
S_{\Gamma}$ and $\delta \in \S_{\Delta},$ we denote by $(\delta,
\gamma, \lambda)$ the permutation $(\lambda(1), \dots,
\lambda(\ell), \gamma(1), \dots, \gamma(i), \delta(1), \dots,
\delta(d-\ell-i)).$ For fixed $\lambda$ and $\delta,$ we denote by
$\tilde{\S}_{\lambda, d, \delta}$ the set of all of the permutations
in the form of $(\lambda, \gamma, \delta).$

\itm In particular, when $\Delta$ is the empty set, i.e., $(\Lambda,
\Gamma)$ is a partition of $[d],$ we simply write
$\tilde{\S}_{\lambda, d, \delta}$ as $\tilde{\S}_{\lambda, d}$ which
is the set of all of permutations in the form of $(\lambda,
\gamma),$ for some fixed $\lambda \in S_{\Lambda}.$

\itm We analogously define $\tilde{\S}_{ d, \delta}$ in the case
that $\Lambda$ is the empty set, i.e., $(\Gamma, \Delta)$ is a
partition of $[d].$
\end{alist}
\end{defn}

\subsection{Right side of (\ref{gdsigma})} Because $\nHat{z}(\sigma,j) = (-1)^j z(\sigma,j)$ and
$\det(\nHat{X}(1, d)) = (-1)^d \det({X}(1, d))$ to prove that the
right side of (\ref{gdsigma}) is equal to $\frac{1}{d!} \sum_{\sigma
\in \S_d} \sgn(\sigma) \prod_{j=1}^d z(\sigma, j)$ is equivalent to
showing that
\begin{equation}\label{det2}
\sum_{\sigma \in \S_d} \sgn(\sigma) \prod_{j=1}^{d}
\nHat{z}(\sigma,j)  = (-1)^{\frac{d(d-1)}{2}} \det(\nHat{X}(\1, d)).
\end{equation}

The following lemma gives a stronger statement. It involves
$\tilde{\S}_{d,\delta}.$ For any $\sigma = (\gamma, \delta) \in
\tilde{\S}_{d, \delta},$ and $\forall 1 \le j \le i,$
$\det(\nHat{X}(\sigma,j))$ and $\det(\nHat{Y}(\sigma,j))$ do not
depend on $\delta.$ So we simply write them as $\det(\nHat{X}(\gamma
,j))$ and $\det(\nHat{Y}(\gamma,j)).$

\begin{lem}\label{ldet} For any  $1 \le i \le d,$ let
$(\Gamma, \Delta)$ be a partition of $[d]$ with the size of $\Gamma$
equal to $i.$ For any $\delta \in \S_{\Delta}$ and $\gamma \in
\S_{\Gamma},$ we have that
\begin{equation}\label{det1}
\sum_{\sigma \in \tilde{\S}_{d, \delta}} \sgn(\sigma)
\prod_{j=1}^{i} \nHat{z}(\sigma,j)  = (-1)^{\frac{i(i-1)}{2}}
\sgn((\gamma, \delta)) \det(\nHat{X}(\gamma,i)).
\end{equation}

In particular, (\ref{det2}) holds.
\end{lem}

\begin{proof}
We prove (\ref{det1}) by induction on $i.$

When $i = 1,$ there is only one $\sigma$ in $\tilde{\S}_{d, \delta}$
and $\sgn(\sigma) = \sgn((\gamma, \delta)).$ Together with the fact
that $\det(\nHat{Y}(\gamma,1)) = 1,$ (\ref{det1}) holds.

Assuming that (\ref{det1}) holds when $i = i_0 \ge 1,$ consider $i =
i_0 +1.$

For any $m: 1 \le m \le i,$ let $\Gamma^{(m)} = \Gamma \setminus \{
\gamma(m)\}$ and $\Delta^{(m)} = \Delta \cup  \{ \gamma(m)\}.$ Then
$(\Gamma^{(m)}, \Delta^{(m)})$ is a partition of $[d],$ where the
size of $\Gamma^{(m)}$ is $i-1 = i_0.$ Let $\gamma^{(m)} =
(\gamma(1), \dots, \gamma(m-1), \gamma(m+1), \dots, \gamma(i))$ and
$\delta^{(m)} = (\gamma(m), \delta(1), \dots, \delta(d-i)).$ We know
that $\sgn((\gamma^{(m)}, \delta^{(m)})) = (-1)^{i+m} \sgn((\gamma,
\delta)).$ Then by the induction hypothesis,
\begin{eqnarray*}
\sum_{\sigma \in \tilde{\S}_{d, \delta^{(m)}}} \sgn(\sigma)
\prod_{j=1}^{i-1} \nHat{z}(\sigma,j)  &=&
(-1)^{\frac{(i-1)(i-2)}{2}} \sgn((\gamma^{(m)}, \delta^{(m)}))
\det(\nHat{X}(\gamma^{(m)},i-1)) \\
&=& (-1)^{\frac{(i-1)(i-2)}{2}+i+m} \sgn((\gamma, \delta))
\det(\nHat{X}(\gamma^{(m)},i-1))
\end{eqnarray*}

However, $(\tilde{\S}_{d, \delta^{(m)}})_{1\le m \le i}$ gives a
partition for $\tilde{\S}_{d, \delta},$ and for any $\sigma \in
\tilde{\S}_{d, \delta},$ $\nHat{z}(\sigma, i)$ is an invariant. In
particular, $\nHat{z}(\sigma, i) = \nHat{z}((\gamma, \delta), i) =
\det(\nHat{X}(\gamma,i))/\det(\nHat{Y}(\gamma,i)).$ Therefore,

\begin{eqnarray*}
& &\sum_{\sigma \in \tilde{\S}_{d, \delta}} \sgn(\sigma)
\prod_{j=1}^{i} \nHat{z}(\sigma,j) =  \sum_{m=1}^i \nHat{z}((\gamma,
\delta), i) \sum_{\sigma \in \tilde{\S}_{d, \delta^{(j)}}}
\sgn(\sigma) \prod_{j=1}^{i-1}
\nHat{z}(\sigma,j)\\
&=& (-1)^{\frac{(i-1)(i-2)}{2}+i-1} \sgn((\gamma,
\delta))\nHat{z}((\gamma, \delta), i) \sum_{m=1}^i  (-1)^{m+1}
\det(\nHat{X}(\gamma^{(m)},i-1)) \\
&=& (-1)^{\frac{i(i-1)}{2}} \sgn((\gamma, \delta))\nHat{z}((\gamma,
\delta), i) \det(\nHat{Y}(\gamma,i)) \\
&=& (-1)^{\frac{i(i-1)}{2}} \sgn((\gamma, \delta))
\det(\nHat{X}(\gamma,i)).
\end{eqnarray*}

Therefore, (\ref{det1}) holds. If we set $i = d,$ then $\Delta =
\emptyset,$ and $\Gamma = [d].$ Letting $\gamma = \1$ be the
identity in $\S_d,$ we obtain (\ref{det2}).
\end{proof}

\subsection{Left side of (\ref{gdsigma})} The proof that
\begin{equation}\label{gsigma3}
\sum_{\sigma \in \S_d} \sgn(\sigma) g_d(z(\sigma,1), \dots,
z(\sigma,d))  = \frac{1}{d!} \sum_{\sigma \in \S_d} \sgn(\sigma)
\prod_{j=1}^d z(\sigma, j)
\end{equation} is relatively harder than what we did in the previous section. We need use
the following lemma.

\begin{lem}\label{zerofive}
For any $0 \le \ell +  k \le d-2,$  given $p(y_1, \dots, y_{\ell})$
a function on $\ell$ variables, let $q(\sigma) = p(z(\sigma,1),
\dots, z(\sigma,\ell)), \forall \sigma \in \S_d.$ Then

\begin{equation}\label{zero5}
\sum_{\sigma \in \S_d} \sgn(\sigma) q(\sigma) \frac{\prod_{j=\ell +
1}^{d} {z}(\sigma,j) }{({z}(\sigma,\ell+1))^{k+1}} = 0.
\end{equation}

\end{lem}

Given this lemma, we are able to prove the following proposition
which implies (\ref{gsigma3}) when we set $\ell = 0.$

\begin{prop} Define $s_0 = 1, z(\sigma, 0)=1.$ For any $\ell: 0 \le \ell \le
d,$ we have that

\begin{eqnarray}\label{gsigma2}
& &\sum_{\sigma \in \S_d} \sgn(\sigma) g_d(z(\sigma,1), \dots,
z(\sigma,d))  \\
&=& \frac{1}{(d-\ell)!} \sum_{\sigma \in \S_d} \sgn(\sigma)
\sum_{s_1=1}^{\frac{z(\sigma,1)}{z(\sigma,0)} s_0} \cdots
\sum_{s_\ell=1}^{\frac{z(\sigma,\ell)}{z(\sigma,\ell-1)} s_{\ell-1}}
\frac{\prod_{j=\ell+1}^d z(\sigma, j)
}{(z(\sigma,\ell))^{d-\ell}}s_\ell^{d-\ell}. \nonumber
\end{eqnarray}
\end{prop}

\begin{proof}
We proceed by descending induction on $\ell.$

When $\ell = d,$ (\ref{gsigma2}) holds by the definition of $g_d.$

When $\ell = d-1,$ it's easy to check that (\ref{gsigma2}) holds.

Assuming (\ref{gsigma2}) holds for $\ell = \ell_0 + 1 \le d-1,$ we
consider $\ell = \ell_0 \le d-2.$ By the induction hypothesis,
\begin{eqnarray*}
& &\sum_{\sigma \in \S_d} \sgn(\sigma) g_d(z(\sigma,1), \dots,
z(\sigma,d))  \\
&=& \frac{1}{(d-\ell-1)!} \sum_{\sigma \in \S_d} \sgn(\sigma)
\sum_{s_1=1}^{\frac{z(\sigma,1)}{z(\sigma,0)} s_0} \cdots
\sum_{s_{\ell+1}=1}^{\frac{z(\sigma,\ell+1)}{z(\sigma,\ell)}
s_{\ell}} \frac{\prod_{j=\ell+2}^d z(\sigma, j)
}{(z(\sigma,\ell+1))^{d-\ell-1}}s_{\ell+1}^{d-\ell-1} \\
&=& \frac{1}{(d-\ell-1)!} \sum_{\sigma \in \S_d} \sgn(\sigma)
\sum_{s_1=1}^{\frac{z(\sigma,1)}{z(\sigma,0)} s_0} \cdots
\sum_{s_\ell=1}^{\frac{z(\sigma,\ell)}{z(\sigma,\ell-1)} s_{\ell-1}}
\frac{\prod_{j=\ell+2}^d z(\sigma, j)
}{(z(\sigma,\ell+1))^{d-\ell-1}}
P_{d-\ell-1}(\frac{z(\sigma,\ell+1)}{z(\sigma,\ell)} s_{\ell}).
\end{eqnarray*}
Recall that $P_{d-\ell-1}(x)$ is the power sum polynomial. Note that
$d-\ell-1 \ge 1.$ By (\ref{pdeg}), (\ref{pconst}) and (\ref{plead}),
we can assume $$P_{d-\ell-1}(x) = \frac{1}{d-\ell}x^{d-\ell} +
\sum_{m = 1}^{d-\ell-1} c_m x^m,$$ where $c_m \in \R.$

For $\forall m: 1 \le m \le d-\ell-1,$ defining $x_0 = 1,$ let
$$p_m(x_1, \dots, x_{\ell}) = \sum_{s_1=1}^{\frac{x_1}{x_0} s_0} \cdots
\sum_{s_\ell=1}^{\frac{x_\ell}{x_{\ell-1}} s_{\ell-1}}
\left(\frac{s_\ell}{x_\ell}\right)^m.$$ Then $p_m$ is a function on
$\ell$ variables. Let
$$q_m(\sigma) = p_m(z(\sigma,1), \dots, z(\sigma,\ell)).$$
Then
\begin{eqnarray*}
& & \sum_{\sigma \in \S_d} \sgn(\sigma)
\sum_{s_1=1}^{\frac{z(\sigma,1)}{z(\sigma,0)} s_0} \cdots
\sum_{s_\ell=1}^{\frac{z(\sigma,\ell)}{z(\sigma,\ell-1)} s_{\ell-1}}
\frac{\prod_{j=\ell+2}^d z(\sigma, j)
}{(z(\sigma,\ell+1))^{d-\ell-1}}
\left(\frac{z(\sigma,\ell+1)}{z(\sigma,\ell)} s_{\ell}\right)^m \\
&=& \sum_{\sigma \in \S_d} \sgn(\sigma) q_m(\sigma)
\frac{\prod_{j=\ell+2}^d z(\sigma, j) }{(z(\sigma,\ell+1))^{d-\ell-1-m}} \\
&=& \sum_{\sigma \in \S_d} \sgn(\sigma) q_m(\sigma)
\frac{\prod_{j=\ell+1}^d z(\sigma, j)
}{(z(\sigma,\ell+1))^{d-\ell-m}} = 0.
\end{eqnarray*}
The last equality is by (\ref{zero5}). Therefore,
\begin{eqnarray*}
& &\sum_{\sigma \in \S_d} \sgn(\sigma) g_d(z(\sigma,1), \dots,
z(\sigma,d))  \\
&=& \frac{1}{(d-\ell-1)!} \sum_{\sigma \in \S_d} \sgn(\sigma)
\sum_{s_1=1}^{\frac{z(\sigma,1)}{z(\sigma,0)} s_0} \cdots
\sum_{s_\ell=1}^{\frac{z(\sigma,\ell)}{z(\sigma,\ell-1)} s_{\ell-1}}
\frac{\prod_{j=\ell+2}^d z(\sigma, j)
}{(z(\sigma,\ell+1))^{d-\ell-1}} \frac{1}{d-\ell}
\left(\frac{z(\sigma,\ell+1)}{z(\sigma,\ell)}
s_{\ell}\right)^{d-\ell} \\
&=& \frac{1}{(d-\ell)!} \sum_{\sigma \in \S_d} \sgn(\sigma)
\sum_{s_1=1}^{\frac{z(\sigma,1)}{z(\sigma,0)} s_0} \cdots
\sum_{s_\ell=1}^{\frac{z(\sigma,\ell)}{z(\sigma,\ell-1)} s_{\ell-1}}
\frac{\prod_{j=\ell+1}^d z(\sigma, j) }{(z(\sigma,\ell))^{d-\ell}}
s_{\ell}^{d-\ell}.
\end{eqnarray*}
\end{proof}

Now we have everything we need to prove Proposition \ref{gdsigma}.
\subsection*{Proof of Proposition \ref{gdsigma}} The proposition follows from (\ref{gsigma3}), (\ref{det2})
and the facts that $\nHat{z}(\sigma,j) = (-1)^j z(\sigma,j)$ and
$\det(\nHat{X}(1, d)) = (-1)^d \det({X}(1, d)).$

%
%
%

\subsection{Proof of Lemma \ref{zerofive}} It remains to prove Lemma \ref{zerofive}, which is
most complicated part of this section. We will break the proof into
several steps. The first lemma we need involves symmetric
polynomials.

A {\it symmetric polynomial} on $d$ variables $y_1, ..., y_d$ is a
polynomial that is unchanged by any permutation of its variables.

\begin{lem}
For any $k \ge 0$, there exist symmetric polynomials
$\phi_{i}^k(y_1, y_2,\dots,y_i)$ on variables $y_1, y_2, \dots, y_i$
for any $i: 1 \le i \le k+2$ and symmetric polynomials
$\varphi_{j}^k(y_1, y_2, \dots, y_j)$ on variables $y_1, y_2, \dots,
y_j$ for any $j: 2 \le j \le k+2,$ so that
\begin{eqnarray}
& &\phi_1^k = 1, \quad \phi_{k+1}^{k} = 1, \quad \phi_{k+2}^k = 0, \label{sym1} \\
& &\varphi_2^k(y_1,y_2) = \sum_{i=0}^k y_1^i y_2^{k-i},
\quad \varphi_{k+2}^k = 1, \label{sym2} \\
 \forall 1 \le i \le k+1, & &\phi_{i}^k(y_1,\dots, y_i)
y_{i+1}^{k+2-i} - \varphi_{i+1}^{k}(y_1, \dots, y_{i+1}) y_{i+1} \label{sym3} \\
&=& -\phi_{i+1}^k(y_1, \dots, y_{i+1})y_1 \cdots y_{i+1}. \nonumber
\end{eqnarray}
\end{lem}

\begin{proof}
Proof by induction on $k.$

When $k=0,$ $\phi_1^0 = 1,  \phi_2^0 = 0, \varphi_2^0 = 1
\Rightarrow \phi_1^0 y_2 - \varphi_2^0 y_2 = - \phi_2^0 y_1 y_2.$

Assume that (\ref{sym1}),(\ref{sym2}) and (\ref{sym3}) hold for $k =
k_0 \ge 0.$ We check the case $k = k_0+1.$

We set $$\phi_i^k(y_1,\dots, y_i) =\begin{cases}1, & \mbox{if
$i = 1,$} \\
\varphi_i^{k-1}(y_1,\dots, y_i), &\mbox{if $2 \le i \le k+1,$}\\
0, & \mbox{if $i = k+2.$}
\end{cases}$$
Note that $\phi_{k+1}^k = \varphi_{k+1}^{k-1} = 1$ by the induction
hypothesis. Thus, (\ref{sym1}) holds.

Now all of the $\phi_{i}^k$'s are given. In order to satisfy
(\ref{sym3}), for $\forall 1 \le i \le k+1$, we set
\begin{equation}\label{varphi}\varphi_{i+1}^{k}(y_1, \dots, y_{i+1})
= \phi_{i}^k(y_1,\dots, y_i) y_{i+1}^{k+1-i} + \phi_{i+1}^k(y_1,
\dots, y_{i+1})y_1 \cdots y_{i}.\end{equation}

Hence, it is left to show that $\varphi_{i+1}^k$'s are symmetric
polynomials and satisfy (\ref{sym2}).

When $i = 1,$ $\varphi_{2}^{k}(y_1, y_{2}) = \phi_{1}^k(y_1)
y_{2}^{k} + \phi_{2}^k(y_1, y_{2})y_1 = y_2^k +
\varphi_{2}^{k-1}(y_1, y_2)y_1 = y_2^k + \left(\sum_{i=0}^{k-1}y_1^i
y_2^{k-1-i}\right) y_1 = \sum_{i=0}^k y_1^i y_2^{k-i}.$

When $2 \le i \le k,$ because the right hand side of (\ref{varphi})
is symmetric on $y_1, y_2, \dots, y_i,$ it's enough to show that it
is symmetric on $y_1$ and $y_{i+1}.$ However,
\begin{eqnarray*}
\phi_{i}^k(y_1,\dots, y_i) &=& \varphi_{i}^{k-1}(y_1, \dots, y_i) \\
&=& \phi_{i-1}^{k-1}(y_1,\dots, y_{i-1}) y_{i}^{k+1-i} +
\phi_{i}^{k-1}(y_1, \dots, y_{i})y_1 \cdots y_{i-1}.
\end{eqnarray*}
Because $\phi_i^k$ is symmetric, we can switch $y_1$ and $y_i.$ So
$$\phi_{i}^k(y_1,\dots, y_i) = \phi_{i-1}^{k-1}(y_2,\dots, y_{i}) y_{1}^{k+1-i} +
\phi_{i}^{k-1}(y_1, \dots, y_{i})y_2 \cdots y_{i}.$$ Similarly,
$$\phi_{i+1}^k(y_1, \dots, y_{i+1}) = \phi_{i}^{k-1}(y_2,\dots, y_{i+1}) y_{1}^{k-i} +
\phi_{i+1}^{k-1}(y_1, \dots, y_{i+1})y_2 \cdots y_{i+1}.$$
Therefore, \begin{eqnarray*} & &\varphi_{i+1}^{k}(y_1, \dots,
y_{i+1}) \\
&=& \phi_{i-1}^{k-1}(y_2,\dots, y_{i}) y_{1}^{k+1-i}y_{i+1}^{k+1-i}
+ \phi_{i}^{k-1}(y_1, \dots, y_{i})y_2 \cdots y_{i}y_{i+1}^{k+1-i} +
\\
& &\phi_{i}^{k-1}(y_2,\dots, y_{i+1}) y_{1}^{k+1-i}y_2 \cdots y_{i}
+ \phi_{i+1}^{k-1}(y_1, \dots, y_{i+1})y_1 y_2^2 \cdots y_{i}^2
y_{i+1}
\end{eqnarray*}
is symmetric on $y_1$ and $y_{i+1}.$

When $i = k+1,$ $\varphi_{k+2}^{k}(y_1, \dots, y_{i+1}) =
\phi_{k+1}^k + \phi_{k+2}^k y_1 \cdots y_{k+1} = 1.$
\end{proof}

\begin{lem}\label{lzero} For any $0 \le k \le d-2,$ $1 \le i \le k+2,$ let
$(\Gamma, \Delta)$ be a partition of $[d]$ with the size of $\Gamma$
equal to $i.$ For any $\delta \in \S_{\Delta}$ and $\gamma \in
\S_{\Gamma},$ we have that
\begin{equation}\label{zero1}
\sum_{\sigma \in \tilde{\S}_{d,\delta}} \sgn(\sigma)
\frac{\prod_{j=1}^{i-1} \nHat{z}(\sigma,j)
}{(\nHat{z}(\sigma,1))^{k+1}} = (-1)^{\frac{i(i+1)}{2}-1}
\sgn((\gamma, \delta)) \frac{\phi_{i}^k( \nHat{x}_{\gamma(1),1},
\dots, \nHat{x}_{\gamma(i),1} ) }{\prod_{j=1}^i
(\nHat{x}_{\gamma(j),1})^{k+2-i}} \det(\nHat{Y}(\gamma,i)).
\end{equation}
\end{lem}

\begin{proof}

We prove (\ref{zero1}) by induction on $i.$

When $i = 1,$ there is only one $\sigma$ in $\tilde{\S}_{d, \delta}$
and $\sgn(\sigma) = \sgn((\gamma, \delta)).$ Together with the facts
that $\phi_1^k = 1,$ $\nHat{z}(\sigma, 1) = \nHat{x}_{\sigma(1),1} =
\nHat{x}_{\gamma(1),1}$ and $\det(\nHat{Y}(\gamma,1)) = 1,$ we
conclude (\ref{zero1}).

Assuming that (\ref{zero1}) holds when $i = i_0 \ge 1,$ consider $i
= i_0 +1.$

For any $m: 1 \le m \le i,$ let $\Gamma^{(m)} = \Gamma \setminus \{
\gamma(m)\}$ and $\Delta^{(m)} = \Delta \cup  \{ \gamma(m)\}.$ Then
$(\Gamma^{(m)}, \Delta^{(m)})$ is a partition of $[d],$ where the
size of $\Gamma^{(m)}$ is $i-1 = i_0.$ Let $\gamma^{(m)} =
(\gamma(1), \dots, \gamma(m-1), \gamma(m+1), \dots, \gamma(i))$ and
$\delta^{(m)} = (\gamma(m), \delta(1), \dots, \delta(d-i)).$ Then by
the induction hypothesis, we have that
\begin{eqnarray*}
& &\sum_{\sigma \in \tilde{\S}_{d, \delta^{(m)}}} \sgn(\sigma)
\frac{\prod_{j=1}^{i-2} \nHat{z}(\sigma,j)
}{(\nHat{z}(\sigma,1))^{k+1}} \\
&=& (-1)^{\frac{i(i-1)}{2}-1} \sgn((\gamma^{(m)}, \delta^{(m)}))
\frac{\phi_{i-1}^k( \nHat{x}_{\gamma^{(m)}(1),1}, \dots,
\nHat{x}_{\gamma^{(m)}(i-1),1} ) }{\prod_{j=1}^{i-1}
(\nHat{x}_{\gamma^{(m)}(j),1})^{k+3-i}}
\det(\nHat{Y}(\gamma^{(m)},i-1)).
\end{eqnarray*}

However, $(\tilde{\S}_{d, \delta^{(m)}})_{1\le m \le i}$ gives a
partition for $\tilde{\S}_{d,\delta}.$ Therefore,

$$\sum_{\sigma \in \tilde{\S}_{d, \delta}} \sgn(\sigma) \frac{\prod_{j=1}^{i-1}
\nHat{z}(\sigma,j) }{(\nHat{z}(\sigma,1))^{k+1}}  = \sum_{m=1}^i
\sum_{\sigma \in \tilde{\S}_{d,\delta^{(j)}}} \sgn(\sigma)
\frac{\prod_{j=1}^{i-2} \nHat{z}(\sigma,j)
}{(\nHat{z}(\sigma,1))^{k+1}} \nHat{z}(\sigma, i-1).$$

But for any $\sigma \in \tilde{\S}_{d, \delta^{(m)}},$
$\nHat{z}(\sigma, i-1)$ is an invariant.  In particular,
$\nHat{z}(\sigma, i-1) = \nHat{z}((\gamma^{(m)}, \delta^{(m)}), i-1)
=
\det(\nHat{X}(\gamma^{(m)},i-1))/\det(\nHat{Y}(\gamma^{(m)},i-1)).$

Hence,
\begin{eqnarray*}
& &\sum_{\sigma \in \tilde{\S}_{d, \delta}} \sgn(\sigma)
\frac{\prod_{j=1}^{i-1}
\nHat{z}(\sigma,j) }{(\nHat{z}(\sigma,1))^{k+1}} \\
&=& \sum_{m=1}^i (-1)^{\frac{i(i-1)}{2}-1} \sgn((\gamma^{(m)},
\delta^{(m)})) \frac{\phi_{i-1}^k( \nHat{x}_{\gamma^{(m)}(1),1},
\dots, \nHat{x}_{\gamma^{(m)}(i-1),1} ) }{\prod_{j=1}^{i-1}
(\nHat{x}_{\gamma^{(m)}(j),1})^{k+3-i}}
\det(\nHat{X}(\gamma^{(m)},i-1)).
\end{eqnarray*}
Note that $(\nHat{x}_{\gamma(m), 1})^{k+3-i} \prod_{j=1}^{i-1}
(\nHat{x}_{\gamma^{(m)}(j),1})^{k+3-i} = \prod_{j=1}^{i}
(\nHat{x}_{\gamma(j),1})^{k+3-i} ,$ and $\sgn((\gamma^{(m)},
\delta^{(m)})) = (-1)^{i+m} \sgn((\gamma, \delta)).$ Therefore,
$$\sum_{\sigma \in \tilde{\S}_{d, \delta}} \sgn(\sigma) \frac{\prod_{j=1}^{i-1}
\nHat{z}(\sigma,j) }{(\nHat{z}(\sigma,1))^{k+1}} \\
= (-1)^{\frac{i(i-1)}{2}-1+i} \frac{\sgn((\gamma,
\delta))}{\prod_{j=1}^{i} (\nHat{x}_{\gamma(j),1})^{k+3-i}} A,$$
where
\begin{eqnarray*}
A &=& \sum_{m=1}^i (-1)^m \phi_{i-1}^k(
\nHat{x}_{\gamma^{(m)}(1),1}, \dots, \nHat{x}_{\gamma^{(m)}(i-1),1}
)
(\nHat{x}_{\gamma(m),1})^{k+3-i} \det(\nHat{X}(\gamma^{(m)},i-1)) \\
&=& - \det\left(\left(\begin{array}{cccccc}
\phi_{i-1}^k(\nHat{x}_{\gamma^{(1)}(1),1}, \dots,
\nHat{x}_{\gamma^{(1)}(i-1),1}) (\nHat{x}_{\gamma(1),1})^{k+3-i} & \nHat{x}_{\gamma(1), 1} &  \cdots & \nHat{x}_{\gamma(1), i-1} \\
\phi_{i-1}^k( \nHat{x}_{\gamma^{(2)}(1),1}, \dots,
\nHat{x}_{\gamma^{(2)}(i-1),1} )
(\nHat{x}_{\gamma(2),1})^{k+3-i} & \nHat{x}_{\gamma(2), 1} &  \cdots & \nHat{x}_{\gamma(2), i-1} \\
\vdots & \vdots &  \ddots & \vdots \\
\phi_{i-1}^k( \nHat{x}_{\gamma^{(i)}(1),1}, \dots,
\nHat{x}_{\gamma^{(i)}(i-1),1} ) (\nHat{x}_{\gamma(i),1})^{k+3-i} &
\nHat{x}_{\gamma(i), 1} &  \cdots & \nHat{x}_{\gamma(i), i-1}
\end{array}\right)\right)
\end{eqnarray*}

By (\ref{sym3}), if we subtract the second column times
$\varphi_{i}^{k}(\nHat{x}_{\gamma(1),1}, \dots,
\nHat{x}_{\gamma(i),1})$ from the first column, then
\begin{eqnarray*}
A &=& - \det\left(\left(\begin{array}{cccccc}
-\phi_{i}^k(\nHat{x}_{\gamma(1),1}, \dots,
\nHat{x}_{\gamma(i),1})\nHat{x}_{\gamma(1),1} \cdots
\nHat{x}_{\gamma(i),1} & \nHat{x}_{\gamma(1), 1} &  \cdots & \nHat{x}_{\gamma(1), i-1} \\
-\phi_{i}^k(\nHat{x}_{\gamma(1),1}, \dots,
\nHat{x}_{\gamma(i),1})\nHat{x}_{\gamma(1),1} \cdots
\nHat{x}_{\gamma(i),1}  & \nHat{x}_{\gamma(2), 1} &  \cdots & \nHat{x}_{\gamma(2), i-1} \\
\vdots & \vdots  & \ddots & \vdots \\
-\phi_{i}^k(\nHat{x}_{\gamma(1),1}, \dots,
\nHat{x}_{\gamma(i),1})\nHat{x}_{\gamma(1),1} \cdots
\nHat{x}_{\gamma(i),1}  & \nHat{x}_{\gamma(i), 1} & \cdots &
\nHat{x}_{\gamma(i), i-1}
\end{array}\right)\right) \\
&=& \phi_{i}^k(\nHat{x}_{\gamma(1),1}, \dots,
\nHat{x}_{\gamma(i),1})\nHat{x}_{\gamma(1),1} \cdots
\nHat{x}_{\gamma(i),1} \det( \nHat{Y}(\gamma, i)).
\end{eqnarray*}

Therefore,

$$\sum_{\sigma \in \tilde{\S}_{d, \delta}} \sgn(\sigma) \frac{\prod_{j=1}^{i-1}
\nHat{z}(\sigma,j) }{(\nHat{z}(\sigma,1))^{k+1}} \\
= (-1)^{\frac{i(i+1)}{2}-1} \sgn((\gamma, \delta))
\frac{\phi_{i}^k(\nHat{x}_{\gamma(1),1}, \dots,
\nHat{x}_{\gamma(i),1})}{\prod_{j=1}^{i}
(\nHat{x}_{\gamma(j),1})^{k+2-i}} \det (\nHat{Y}(\gamma, i)).$$
\end{proof}

\begin{cor}
For any $0 \le \ell +  k \le d-2,$ $ i = k+2,$ let $(\Lambda,
\Gamma, \Delta)$ be a partition of $[d]$ with the sizes of $\lambda$
and $\Gamma$ equal to $\ell$ and $i,$ respectively. For any $\delta
\in \S_{\Delta}$ and any $\lambda \in \S_{\Lambda},$  we have that
\begin{equation}\label{zero2}
\sum_{\sigma \in \tilde{\S}_{\lambda, d, \delta}} \sgn(\sigma)
\frac{\prod_{j=\ell + 1}^{\ell+k+1} \nHat{z}(\sigma,j)
}{(\nHat{z}(\sigma,\ell+1))^{k+1}} = 0,
\end{equation}
and
\begin{equation}\label{zero3}
\sum_{\sigma \in \tilde{\S}_{\lambda, d, \delta}} \sgn(\sigma)
\frac{\prod_{j=\ell + 1}^{\ell+k+1} {z}(\sigma,j)
}{({z}(\sigma,\ell+1))^{k+1}} = 0,
\end{equation}
\end{cor}

\begin{proof}

Proof by induction on $\ell.$

When $\ell = 0,$ (\ref{zero2}) follows from Lemma \ref{lzero} and
the fact that $\phi_{k+2}^k = 0.$ Thus, (\ref{zero3}) holds by
(\ref{zhat}).

We assume for $\ell = \ell_0 \ge 0,$ (\ref{zero2}) and (\ref{zero3})
holds. We check the case $\ell = \ell_0 +1.$ Because (\ref{zero2})
and (\ref{zero3}) are equivalent by (\ref{zhat}), it's enough to
show (\ref{zero2}).

Without loss of generality, we assume that $d \in \Lambda$ and
$\lambda(1) = d.$ For $1 \le q \le d-1,$ define
$$y_{p,q} = \begin{cases} {(\nHat{x}_{p,q+1}-
\nHat{x}_{d,q+1})}/{(\nHat{x}_{p,1}- \nHat{x}_{d,1})}, & \mbox{if $1
\le p \le d-1$} \\
{\nHat{x}_{p,q+1}}/{\nHat{x}_{p,1}}, & \mbox{if $p = d$}.
\end{cases}$$
Let $W$ be the vertex set $ \{ w_1, w_2, \dots, w_{d-1} \},$ where
the coordinates of $w_p$ are $(y_{p,1}, y_{p,2}, \dots, y_{p,d-1}).$
For any $\sigma \in \tilde{\S}_{\lambda, \delta},$ let $\varsigma =
(\sigma(2), \sigma(3), \dots, \sigma(d)).$ Because $\lambda(1) = d,$
$\varsigma \in \S_{d-1}.$ Clearly, $(\Lambda \setminus \{ d \},
\Gamma, \Delta)$ is a partition for $[d-1]$ and $\varsigma \in
\tilde{\S}_{\lambda', d-1, \gamma},$ where $\lambda' = (\lambda(2),
\dots, \lambda(\ell)).$ Therefore, for $j \ge 2,$

\begin{eqnarray*}
& &\det(\nHat{X}(\sigma, j)) \\
&=& (-1)^{j-1} \det \left(\left(\begin{array}{ccccc}
\nHat{x}_{\sigma(2), 1} & \nHat{x}_{\sigma(2),2} & \cdots & \nHat{x}_{\sigma(2), j} \\
 \vdots & \vdots & \ddots & \vdots \\
\nHat{x}_{\sigma(j), 1} & \nHat{x}_{\sigma(j),2} & \cdots &
\nHat{x}_{\sigma(j), j} \\
\nHat{x}_{d, 1} & \nHat{x}_{d,2} & \cdots & \nHat{x}_{d, j} \\
\end{array}\right)\right) \\
&=& (-1)^{j-1} \det \left(\left(\begin{array}{ccccc}
\nHat{x}_{\sigma(2), 1} - \nHat{x}_{d, 1} & \nHat{x}_{\sigma(2),2} -
\nHat{x}_{d,2}
& \cdots & \nHat{x}_{\sigma(2), j}-\nHat{x}_{d, j} \\
 \vdots & \vdots & \ddots & \vdots \\
\nHat{x}_{\sigma(j), 1} - \nHat{x}_{d, 1} & \nHat{x}_{\sigma(j),2} -
\nHat{x}_{d,2} & \cdots &
\nHat{x}_{\sigma(j), j} -\nHat{x}_{d, j}\\
\nHat{x}_{d, 1} & \nHat{x}_{d,2} & \cdots & \nHat{x}_{d, j} \\
\end{array}\right)\right) \\
&=& (-1)^{j-1}\nHat{x}_{d, 1} \prod_{p=2}^{j} (\nHat{x}_{\sigma(p),
1} - \nHat{x}_{d, 1}) \det \left(\left(\begin{array}{ccccc} 1 &
y_{\varsigma(1),1} & \cdots & y_{\varsigma(1), j-1} \\
 \vdots & \vdots & \ddots & \vdots \\
1 & y_{\varsigma(j-1),1} & \cdots & y_{\varsigma(j-1),j-1}\\
1 & y_{d-1,2} & \cdots & y_{d-1, j-1} \\
\end{array}\right)\right) \\
&=& (-1)^{j-1}\nHat{x}_{d, 1} \prod_{p=2}^{j} (\nHat{x}_{\sigma(p),
1} - \nHat{x}_{d, 1}) \det(X_W(\varsigma, j-1)).
\end{eqnarray*}
Similarly, $$\det(\nHat{Y}(\sigma, j)) = \prod_{p=2}^{j}
(\nHat{x}_{\sigma(p), 1} - \nHat{x}_{d, 1}) \det(Y_W(\varsigma,
j-1)).$$ Hence,
$$\nHat{z}(\sigma, j) = \det(\nHat{X}(\sigma, j))/\det(\nHat{Y}(\sigma, j))=(-1)^{j-1}\nHat{x}_{d, 1} z_W(\varsigma, j-1).$$

Note that $\sgn(\sigma) = (-1)^{d-1} \sgn(\varsigma).$ Hence, by the
induction hypothesis,
\begin{eqnarray*}
& &\sum_{\sigma \in \tilde{\S}_{\lambda, d, \delta}} \sgn(\sigma)
\frac{\prod_{j=\ell + 1}^{\ell+k+1} \nHat{z}(\sigma,j)
}{(\nHat{z}(\sigma,\ell+1))^{k+1}} \\
&=& \sum_{\varsigma \in \tilde{\S}_{\lambda', d-1, \delta}}
(-1)^{d-1} \sgn(\varsigma) \frac{\prod_{j=\ell + 1}^{\ell+k+1}
(-1)^{j-1} {z}_W(\varsigma,j-1) }{({z}_W(\varsigma,\ell))^{k+1}} =
0.
\end{eqnarray*}
\end{proof}

\begin{proof}[Proof of Lemma \ref{zerofive}]

Consider any partition $(\Lambda, \Gamma, \Delta)$ of $[d],$ where
the size of $\Lambda$ is $\ell$ and the size of $\Gamma$ is $i =
k+2.$ If we fix $\lambda \in \S_{\Lambda}$ and $\delta \in
\S_{\Delta},$ then $\forall \sigma \in \tilde{\S}_{\lambda,d,
\delta},$ $z(\sigma, j)$ is an invariant when $1 \le j \le \ell$ or
$\ell + i \le j \le d.$ Therefore, by (\ref{zero3}),
$$\sum_{\sigma \in \tilde{\S}_{\lambda,d, \delta}} \sgn(\sigma) q(\sigma)
\frac{\prod_{j=\ell + 1}^{d} {z}(\sigma,j)
}{({z}(\sigma,\ell+1))^{k+1}} = 0.$$ But all of the
$\tilde{\S}_{\lambda,d, \delta}$'s give a partition for $\S_d.$
Thus, (\ref{zero5}) holds.
\end{proof}

\section{Examples and Further discussion}

\subsection{Examples of lattice-face polytopes}
In this subsection, we use a fixed family of lattice-face polytopes
to illustrate our results. Let $d = 3,$ and for any positive integer
$k,$ let $P_k$ be the polytope with the vertex set $V =
\{v_1=(0,0,0), v_2=(4,0,0), v_3=(3,6,0), v_4=(2,2,10k)\}.$ One can
check that $P_k$ is a lattice-face polytope.

\begin{ex}[Example of Theorem \ref{main1}]  The volume of $P_k$ is $40k,$ and
$$i(P_k, m) = 40k m^3 + 12 m^2 + 4 m + 1.$$
$\pi(P_k) = \conv\{(0,0), (4,0), (3,6)\},$ where
$$i(\pi(P_k), m) = 12 m^2 + 4m + 1.$$ So
$$i(P_k, m) = 40k m^3 + i(\pi(P_k), m),$$ which agrees with Theorem
\ref{main1}.
\end{ex}

\begin{ex}[Example of Formula (\ref{esum})]

$F_4 = \conv(v_1,v_2,v_3)$ is a negative facet. The hyperplane
determined by $F_4$ is $H = \{(x_1, x_2, x_3) \ | \ x_3 = 0\}.$
Thus, $v_4' = \pi^{-1}(\pi(v_4))\cap H = (2,2,0).$

$F_3 = \conv(v_1, v_2, v_4 )$ is a positive facet. $\pi(F_3) =
\conv((0,0),(4,0),(2,2)).$ $\Omega(\pi(F_3)) = \pi(F_3) \setminus
\conv((0,0), (4,0)).$ $F_3' = \pi^{-1}(\pi(F_3))\cap H = \conv(v_1,
v_2, v_4' ).$ So $$Q_3 = \conv(F_3 \cup F_3') = \conv(v_1, v_2, v_4,
v_4'),$$  $$\rho^+(\Omega(\pi(F_3)), Q_3) = Q_3 \setminus F_3'.$$

$F_2 = \conv(v_1, v_3, v_4 )$ is a positive facet. $\pi(F_2) =
\conv((0,0),(3,6),(2,2)).$ $\Omega(\pi(F_2)) = \pi(F_2) \setminus
(\conv((0,0), (2,2)) \cup \conv((2,2),(3,6))).$ $F_2' =
\pi^{-1}(\pi(F_2))\cap H = \conv(v_1, v_3, v_4' ).$ So $$Q_2 =
\conv(F_2 \cup F_2') = \conv(v_1, v_3, v_4, v_4'),$$
$$\rho^+(\Omega(\pi(F_2)), Q_2) = Q_2 \setminus (F_2' \cup \conv(v_1,
v_4, v_4')\cup \conv(v_3, v_4, v_4')).$$

$F_1 = \conv(v_2, v_3, v_4 )$ is a positive facet. $\pi(F_1) =
\conv((4,0),(3,6),(2,2)).$ $\Omega(\pi(F_1)) = \pi(F_1) \setminus
\conv((4,0), (2,2)).$ $F_1' = \pi^{-1}(\pi(F_1))\cap H = \conv(v_2,
v_3, v_4' ).$ So $$Q_1 = \conv(F_1 \cup F_1') = \conv(v_2, v_3, v_4,
v_4'),$$ $$\rho^+(\Omega(\pi(F_1)), Q_1) = Q_1 \setminus (F_1' \cup
\conv(v_2, v_4, v_4')).$$

Therefore, $$\Omega(P_k) = P_k \setminus F_4 = - \sgn(F_4)
\bigoplus_{i=1}^{3} \sgn(F_i) \rho^+(\Omega(\pi(F_i)), Q_i),$$ which
agrees with Proposition \ref{psum}.
\end{ex}

\begin{ex}[Example of decomposition] In this example, we decompose
$P_k$ into $3!$ sets, where $5$ of them have positive signs and one
has negative sign, which is different from the cases for cyclic
polytopes, where half of the sets have positive signs and the other
half have negative signs.

 Recall that $v_{\sigma,3} = v_4 =
(2,2,10k),$ for any $\sigma \in \S_3.$

When $\sigma = 123 \in \S_3,$ $v_{123, 2} = v_4' = (2,2,0),$
$v_{123,1} = (2,0,0)$ and $v_{123,0} = v_1 = (0,0,0).$ Then
$$S_{123} = \conv(\{v_{123,i} \}_{0\le i \le 3})  \setminus
\conv(\{v_{123,i} \}_{0\le i \le 2}),$$ with $\sgn(123, P_k) = +1.$

When $\sigma = 213 \in \S_3,$ $v_{213, 2} = v_4' = (2,2,0),$
$v_{213,1} = (2,0,0)$ and $v_{213,0} = v_2 = (4,0,0).$  Then
$$S_{213} = \conv(\{v_{213,i} \}_{0\le i \le 3}) \setminus
(\conv(\{v_{213,i} \}_{0\le i \le 2}) \cup \conv(\{v_{213,i}
\}_{1\le i \le 3})),$$ with $\sgn(213, P_k) = +1.$

One can check that $$S_{123} \oplus S_{213} =
\rho^+(\Omega(\pi(F_3)), Q_3).$$

When $\sigma = 231 \in \S_3,$ $v_{231, 2} = v_4' = (2,2,0),$
$v_{231,1} = (2,12,0)$ and $v_{231,0} = v_2 = (4,0,0).$ Then
$$S_{231} = \conv(\{v_{231,i} \}_{0\le i \le 3})  \setminus
(\conv(\{v_{231,i} \}_{0\le i \le 2}) \cup \conv(\{v_{231,i}
\}_{i=0,2,3} \cup \conv(\{v_{231,i} \}_{1 \le i \le 3})),$$ with
$\sgn(231, P_k) = +1.$

When $\sigma = 321 \in \S_3,$ $v_{321, 2} = v_4' = (2,2,0),$
$v_{321,1} = (2,12,0)$ and $v_{321,0} = v_3 = (3,6,0).$ Then
$$S_{321} = \conv(\{v_{321,i} \}_{0\le i \le 3})  \setminus
(\conv(\{v_{321,i} \}_{0\le i \le 2}) \cup \conv(\{v_{321,i}
\}_{i=0,2,3} \cup \conv(\{v_{321,i} \}_{1 \le i \le 3})),$$ with
$\sgn(321, P_k) = -1.$

One can check that $$S_{231} \ominus S_{321} =
\rho^+(\Omega(\pi(F_1)), Q_1).$$

Similarly, we have that $$S_{132} \oplus S_{312} =
\rho^+(\Omega(\pi(F_2)), Q_2).$$

Therefore, $\Omega(P_k) = \bigoplus_{\sigma \in \S_3} \sgn(\sigma,
P_k) S_{\sigma},$ which coincides with Theorem \ref{decomp}.
\end{ex}
\subsection{Further discussion}
Recall that Remark \ref{dflf2} gives an alternative definition for
lattice-face polytopes. Note that in this definition, when $k=0,$
satisfying (\ref{dfcon}) is equivalent to saying that $P$ is an
integral polytope, which implies that the last coefficient of the
Ehrhart polynomial of $P$ is $1.$ Therefore, one may ask

\begin{ques}If $P$ is a polytope that satisfies (\ref{dfcon}) for
all $ k \in K,$ where $K$ is a fixed subset of $\{0, 1, \dots,
d-1\},$ can we say something about the Ehrhart polynomials of $P$?
\end{ques}

A special set $K$ can be chosen as the set of consecutive integers
from $0$ to $d',$ where $d'$ is an integer no greater than $d-1.$
Based on some examples in this case, the Ehrhart polynomials seems
to follow a certain pattern, so we conjecture the following:

\begin{conj} Given $d' \le d-1,$ if $P$ is a $d$-polytope with
vertex set $V$ such that $\forall k: 0 \le k \le d',$ (\ref{dfcon})
is satisfied,
then for $0 \le k \le d'$, the coefficient of $m^k$ in $i(P,m)$ is
the same as in $i(\pi^{d-d'}(P),m).$ In other words,
$$i(P,m) = i(\pi^{d-d'}(P),m) + \sum_{i=d'+1}^d c_i m^i .$$
\end{conj}

When $d'=0$, the condition on $P$ is simply that it is integral. And
when $d'=d-1,$ we are in the case that $P$ is a lattice-face
polytope. Therefore, for these two cases, this conjecture is true.


\bibliographystyle{amsplain}
\bibliography{gen}

\end{document}